\documentclass{amsart}

\usepackage{amssymb}
\usepackage{amscd}
\usepackage{dsfont}
\usepackage{enumerate, xspace}
\usepackage{hyperref}

\numberwithin{equation}{section}

\newtheorem{thm}{Theorem}[section]
\newtheorem{lemma}[thm]{Lemma}

\newtheorem{prop}[thm]{Proposition}

\newtheorem{conv}[thm]{Convention}

\newtheorem{ass}{Assumption}

\newtheorem{rem}[thm]{Remark}

\newcommand\bC{{\mathbb C}}
\newcommand\bE{{\mathbb E}}
\newcommand\bN{{\mathbb N}}
\newcommand\bP{{\mathbb P}}
\newcommand\bR{{\mathbb R}}

\newcommand\bZ{{\mathbb Z}}

\newcommand{\cB}{{\mathcal B}}
\newcommand{\cC}{{\mathcal C}}

\newcommand{\cF}{{\mathcal F}}

\newcommand{\cL}{{\mathcal L}}

\newcommand{\cM}{{\mathcal M}}
\newcommand{\cN}{{\mathcal N}}
\newcommand{\cO}{{\mathcal O}}

\newcommand{\cT}{{\mathcal T}}
\newcommand{\cX}{{\mathcal X}}

\newcommand\btheta{{\boldsymbol{\theta}}}
\newcommand\bomega{{\boldsymbol{\omega}}} 
\newcommand\bmu{{\boldsymbol{\mu}}} 
\newcommand\bnu{{\boldsymbol{\nu}}} 
\def\btau{{\boldsymbol{\tau}}}

\newcommand{\mue}{{\mu^e}}
\newcommand\Pro{{\operatorname{Pr}\,}}
\newcommand{\Czero}{C_0}
\newcommand{\Cone}{C_1}
\newcommand{\Ctwo}{C_2}

\newcommand{\Cfour}{C_7}
\newcommand{\Cfive}{C_6}
\newcommand{\Csix}{C_4}
\newcommand{\Cseven}{C_3}
\newcommand{\Ceight}{C_5}
\newcommand{\Cnine}{C_8}

\newcommand{\vf}{\varphi}
\newcommand{\ve}{\varepsilon}

\newcommand\Id{\mathds{1}}

\newcommand\ocB{{\overline{\cB}}}
\newcommand\para{{t}}
\newcommand{\tX}{\tilde{X}}

\newcommand{\tg}{\tilde{g}}
\newcommand{\tgw}{{\tilde{g}_w}}
\newcommand{\tDelta}{\tilde{\Delta}}

\newcommand\Var{{\Sigma^2}}
\newcommand\Vartwo{{\Sigma^2_2}}
\newcommand\Card{{\rm Card}}

\def\bv{{\mathbf{v}}}
\def\la{\langle}
\def\ra{\rangle}
\newcommand\txi{{\tilde\xi}}
\newcommand\tcX{{\tilde{\cX}}}
\newcommand\hX{\widehat{X}}
\newcommand\muinv{{\mu^w}}
\newcommand\bmuinv{{\boldsymbol \mu}^w}
\title[Random Walk in deterministic environment]{Random
Walk in deterministically changing environment} 
\author{Dmitry Dolgopyat}
\address{Dmitry Dolgopyat\\
Department of Mathematics\\
University of Maryland\\
4417 Mathematics Bldg,  College Park,  MD 20742, USA}
\email{{\tt dmitry@math.umd.edu}}
\author{Carlangelo Liverani}
\address{Carlangelo Liverani\\
Dipartimento di Matematica\\
II Universit\`{a} di Roma (Tor Vergata)\\
Via della Ricerca Scientifica, 00133 Roma, Italy.}
\email{{\tt liverani@mat.uniroma2.it}}
\thanks{The second author wishes to thank Penn State University where he was visiting during part of this work. We thank the anonymous referees for a very thorough and helpful job.}
\begin{document}
\begin{abstract}
We consider a random walk with transition probabilities weakly dependent on an environment with a deterministic, but strongly chaotic, evolution. We prove that for almost all initial conditions of the environment the walk satisfies the CLT.
\end{abstract}
\maketitle
\section{Introduction}
\label{sec:intro}
The continuing interest in the limit properties of random walks has
generated a remarkable amount of literature (see \cite{Sz, Z} for a
review of the field). In particular, many papers have addressed the
case of random walks in dynamical environments. Apart from few papers
in which special hypotheses are imposed on the form of the transition
probabilities (implying that the process is reversible with respect to
the stationary measure of the environment, see \cite{KV}) the authors
have usually investigated the case in which the evolution of the
environment is described by a Markov process with positive transition probabilities
and the transition probabilities of the walk 
 are close to constant (see, e.g., \cite{BMP, BZ} and references
therein). While such a  situation has recently been settled in great
generality (\cite{DKL}) the case of more complex local dynamics 
and/or far from constant transition probabilities is still wide open.

In this paper we address the first issue establishing conditions under
which the CLT holds for a deterministic local evolution with strong
chaotic properties. Such a CLT is established for almost all the
initial conditions of the environment with respect to a (natural)
stationary measure (this is commonly called a {\em quenched CLT}).  

The deterministic dynamics is taken to be independent at each site
(although our method can be easily extended to weakly interacting cases
(cf. \cite{KL0, DKL})). The single site dynamics is a piecewise expanding
one dimensional map. While multidimensional expanding dynamics could
be treated similarly the case of Anosov map poses
a real problem. Indeed the technique used to control the environment
dynamics is borrowed from the study of coupled map lattices (more
precisely from \cite{KL}) and the extension of such a technique to
coupled Anosov systems is still missing. 
In general, the extension to more general dynamics (with substantially
weaker ergodic properties) would be of interest, but to obtain results
in this direction new ideas seem to be needed. 

Note that the present strategy differs, in its probabilistic part, from the one used in \cite{DKL}. In particular, it is not necesary to prove absolutely continuity of the invariant measure of the environment as seen from the particle with respect to the invariant measure of the environment in a fixed reference frame. We hope that this simplification may be helpful in treating more general cases.

The plan of the paper is the following. In section \ref{sec:model} the
system under investigation is explained in detail and the main result
of the paper is precisely stated. This main result (Theorem \ref{thm:main}) follows after estimating the asymptotic independence of two random walks in the same environment (Lemma \ref{lm:hot}).
In Section \ref{sec:proof_env} we study ergodic properties of the environment. 
Section \ref{sec:invariance} contains the proof of the annealed
(averaged) invariance principle. In section \ref{sec:clt-proof} the
proof of Lemma \ref{lm:hot} is reduced to an estimate on the number of
close encounters (Lemma \ref{lem:NL}). Lemma \ref{lem:NL} 
is proven in Section \ref{sec:twowalks}.
\begin{conv}
\label{conv:contants}
In this paper we will use $C$ to designate a generic constant
depending only on the quantities appearing in the Assumptions \ref{ass:assone}, \ref{ass:two}, \ref{ass:pert}, \ref{ass:ellipticity}
below. We will use instead $C_{a,b,c,\dots}$ for constants depending
also on the parameters $a,b,c,\dots$. Consequently the actual
numerical value of such constants may vary from one occurrence to the
next. On the contrary we will use $C_1, C_2, \dots$, to designate
constants whose value is fixed through the paper.
\end{conv}

\section{Model and Results}
\label{sec:model}
Let $I=[0,1]$ and $T:I\to I$ be a piecewise $\cC^2$ topologically mixing
map such that $|D_xT|\geq \lambda>2$ for each $x\in I$ for which
the derivative is well defined. Then, $I^{\bZ^d}=:\Theta$ is the space of environments (it is a measurable space with the product (Borel) $\sigma$-algebra $\cT$) and $\theta\in\Theta$ is an
environment on $\bZ^d$. This environment evolves deterministically
according to a map $F:\Theta\to\Theta$. 
\begin{ass}[\em Environment Dynamics]\label{ass:assone}
For each $\theta\in\Theta$
\[
(F(\theta))_q:=T(\theta_q).
\]
That is, the evolution is independent at each site.\footnote{The case of weakly coupled maps can be treated similarly by using the techniques introduced in \cite{KL} and used 
here to study the present, simpler, case.}
\end{ass}

The evolution of the environment can be thus seen as a deterministic
Markov process on the space $\Theta^{\bN}=:\Omega$ such that, for all
$(\theta^n)_{n\in\bN}:=\btheta\in\Omega$, $\theta^n:=F^n(\theta^0)$.
If $\mu_0$ is the unique absolutely continuous invariant measure of
$T$,\footnote{For the existence and uniqueness (among measures absolutely continuous w.r.t. Lebesgue) see \cite{Babook}.} then $\mue:=\otimes_{p\in\bZ^d}\mu_0$ is the natural invariant
measure for $F$ we are interested in. In fact, 
it is possible to show (\cite{KL}) that it is the only invariant measure 
in a reasonably large class of measures; see the precise statement below. 

We consider a bounded increment random walk $X_n$ in such an
environment. More precisely, let $\Lambda:=\{ z\in\bZ^d\;:\; \|z\|\leq
\Czero\}$ and $\Delta_n=X_{n+1}-X_n$, then the process is defined
by the transition probabilities
\begin{equation}
\label{eq:walk}
\bP(\{\Delta_n=z\}\;|\;X_n,\theta^0)=\pi_z(\tau^{X_n}\theta^{n}) 
\end{equation}
where $\pi_z\equiv 0$ for all $z\not\in\Lambda$,  $\pi_z(\theta)$ depends only on
$\{\theta_q\}_{q\in\Lambda}$, and, for each $z\in\bZ^d$,
$(\tau^z\theta)_i:=\theta_{i+z}$. We will be interested in the measure
${\bf P}_\nu$ on $\Omega\times (\bZ^d)^\bN$ determined by the above process
when the environment is started with the measure $\nu$ and the walk
starts from zero. We will use the notation ${\bf P}^e$ for ${\bf P}_{\mue}$. Finally, we will use $\bE$ for the expectation with respect to the latter measure and $\bE_\nu$ for the expectations with respect to the process ${\bf P}_\nu$.

\begin{ass}[\em Regularity]\label{ass:two}
The functions $\{\pi_z\}_{z\in\Lambda}$ belong to $\cC^1$.
\end{ass}

The next assumption depends on a parameter $\ve>0$. 
\begin{ass}[\em Perturbative regime]\label{ass:pert} There exists
$\{a_z\}_{z\in\Lambda}\subset \bR_+$, $\sum_{z\in\Lambda}a_z=1$, such that
\[
\|\pi_z-a_z\|_{\cC^1}\leq a_z\ve.
\]
\end{ass}
In the following, when we will say ``assumption \ref{ass:pert} holds for $\ve_i$" we will mean that it holds with $\ve=\ve_i$.
The values $\ve_i$ will be taken small enough for
Theorem \ref{thm:mixing}, Proposition \ref{lem:clt-averaged},  Lemma \ref{lem:up} and Lemma \ref{lem:variation}  to hold.

\begin{ass}[\em Ellipticity]\label{ass:ellipticity} 
For each $l\in\bZ^d\setminus\{0\}$, the function $\left|\sum\limits_{z\in\Lambda}\pi_z e^{i\langle l,z\rangle}\right|\in\cC^0(I^\Lambda, \bR_+)$ is not identically equal to $1$.
\end{ass}

It is well known that to study the properties
of $X_n$ it is convenient to study the process of the environment as
seen from the particle. In fact, such a process can be considered in
several fashions of which the following will be relevant in the
sequel.

\subsection{Process of the environment as seen from the
particle}
Consider the process $\bomega=:(\omega^n)_{n\in\bN}\in\Omega$ described
by the action of the Markov operator $S:L^\infty(\Theta)\to L^\infty(\Theta)$ defined by
\begin{equation}
\label{eq:semigroup-env}
Sf(\omega):=\sum_{z\in\Lambda}\pi_z(\omega)f\circ F(\tau^z\omega)=:\sum_{z\in\Lambda}S_z f.
\end{equation}
\begin{rem}
It is easy to verify that the process $\bomega$, $\omega^0=\theta$, has the same
distribution as the process $(\tau^{X_n}\theta^n)_{n\in\bN}$, $\theta^0=\theta$. 
\end{rem}

We can then consider the measure $\bP_\nu$ on $\Omega$ of the
associated Markov process started with a measure $\nu$.

In analogy with the techniques used in the study of
coupled map lattices \cite{KL} it is then natural to restrict the
space of measures on which $S'$ acts.\footnote{As usual the dual operator $S'$ is defined as $S'\nu(f)=\nu(Sf)$ for all $f\in \cC^0$.} To this end we start
by defining the following norms
\begin{equation}
\label{eq:norms}
\begin{split}
&|\mu|:=\sup_{|\vf|_{\cC^0(\Theta,\bR)}\leq 1}\mu(\vf)\\
&\|\mu\|:=\sup_{i\in\bZ^d}\sup_{|\vf|_{\cC^0(\Theta,\bR)}\leq 1}\mu(\partial_{\theta_i}\vf)
\end{split}
\end{equation}
We then consider the Banach space of complex valued
measures\footnote{By $\cM(\Theta)$ we designate the set of complex
valued finite Borel measures on $\Theta$.} 
\begin{equation}\label{eq:banach}
\cB:=\{\mu\in\cM(\Theta)\;:\; \|\mu\|<\infty\}. 
\end{equation}
It is easy to check that such measures have finite dimensional marginals absolutely
continuous w.r.t. Lebesgue and the densities are functions of bounded
variations with variations bounded by the norm of the measure. Moreover 
$\mue$ is the unique invariant measure for $F$ belonging to $\cB$, \cite{KL}.
\begin{thm}
\label{thm:mixing}
For each dynamics $F$ satisfying assumption \ref{ass:assone} and transition probabilities satisfying assumption \ref{ass:two}, the operator $S'$, is a bounded operator on $\cB$. In addition, there exists $\ve_0>0$, depending on $F$, such that if assumption \ref{ass:pert} holds for $\ve_0$, then there exists a unique invariant probability 
measure $\muinv\in\cB$ ($S'\muinv=\muinv$). This
measure enjoys the following properties: There exists $\eta\in(0,1)$ such that for each $\nu\in\cB$ and local functions $\varphi,\phi\in\cC^0$ each depending only on $L$ variables with the two sets of dependency 
having distance at least $M$
\begin{enumerate}
\item $|\nu(S^n\phi)-\muinv(\phi)\nu(1)|\leq CL\eta^n|\phi|_{\infty}\|\nu\|$
\item $|\muinv(\varphi\phi)-\muinv(\varphi)\muinv(\phi)|\leq CL\eta^{\frac{M}{2\Czero}}|\varphi|_{\infty}|\phi|_\infty$.
\end{enumerate}
\end{thm}
The proof of the above theorem can be found in Section \ref{sec:mixone}.
\begin{rem}\label{rem:stationary}
Theorem \ref{thm:mixing} implies that the process $\bP_{\muinv}$ is a stationary (and ergodic) process.
\end{rem}
\subsection{Annealed statistical properties}\label{subsec:averaged}

\begin{prop}
\label{lem:clt-averaged}
For each dynamics $F$ satisfying assumption \ref{ass:assone} and transition probabilities satisfying assumption \ref{ass:two}, if assumption \ref{ass:pert} is satisfied for $\ve_0$ (where $\ve_0>0$ is as in Theorem \ref{thm:mixing}), then there exists a vector $v\in\bR^d$ and a matrix $\Var\geq0$ such
that, for each probability measure $\nu\in\cB$ we have
\[
\begin{split}
&\frac 1N \bE_\nu(X_N)\to v \\
&\frac{X_N-vN}{\sqrt{N}}\Rightarrow 
\cN\left(0, \Var\right)\quad under \ {\bf P}_\nu.
\end{split}
\] 
Moreover, there exists $\Cone>0$ such that,
setting $\tX_N:=X_N-vN$,  the following 
inequality holds for all $N\in\bN$ and $\para\in\bR^d$:
\[
\left|\bE_\nu\left(e^{\frac i{\sqrt N}\langle \para,\tilde
X_N\rangle}\right)-e^{-\frac12\langle\para,\Var\para\rangle}\right|\leq
\Cone(1+\|\para\|^3)N^{-\frac 12}\|\nu\|.
\]
Finally, if Assumption \ref{ass:ellipticity} is also satisfied, then $\Var>0$.
\end{prop}
\begin{proof}
Let us start noticing that
\[
\frac 1N\bE_\nu\left(X_N\right)=\frac 1N\sum_{k=0}^{N-1}\bE_\nu(\Delta_k)=
\frac 1N\sum_{k=0}^{N-1}\bE_\nu\left(\bE_\nu(\Delta_k\;|\;\cF_k)\right),
\]
where $\cF_k:=\sigma\{\theta^0,X_1,\dots,X_k\}$.
The relevance of the process as seen from the particle is due to the following fact:
\begin{equation}\label{eq:S}
\bE_\nu(\Delta_k\;|\;\cF_k)=\sum_{z\in\Lambda}z\,\pi_z(\tau^{X_{k}}\theta^k)=\sum_{z\in\Lambda}z\,\pi_z(\omega^k)=: g(\omega^k).
\end{equation}
Thus,
\[
\frac 1N\bE_\nu\left(X_N\right)=\frac 1N\sum_{k=0}^{N-1}\nu(g(\omega^k))=\frac 1N\sum_{k=0}^{N-1}[(S')^k\nu]( g).
\]
Accordingly, Theorem \ref{thm:mixing} implies
\begin{equation}\label{eq:average-one}
\lim_{N\to\infty}\frac 1N\bE_\nu\left(X_N\right)=\muinv(g)=:v.
\end{equation}

To prove the CLT let $\tDelta_n=\tX_{n+1}-\tX_n,$ then
\[
\bE_\nu\left(e^{\frac{i}{\sqrt
N}\langle \para,\tilde X_N\rangle}\right)=\bE_\nu\left(e^{\frac{i}{\sqrt
N}\langle \para,\tilde X_{N-1}\rangle}\bE_\nu\left(e^{\frac{i}{\sqrt
N}\langle \para,\Delta_{N-1}-v\rangle}\;\big|\;\cF_{N-1}\right)\right).
\]
Since
\[
\bE_\nu\left(e^{\frac{i}{\sqrt
N}\langle \para,\Delta_{k}-v\rangle}\;\big|\;\cF_{k}\right)=\sum_{z\in\Lambda}\pi_z(\tau^{X_k}\theta^k)e^{\frac{i}{\sqrt
N}\langle \para,z-v\rangle}
\]
it is natural to introduce the operators, for all $t\in\bC^d$,
\begin{equation}\label{eq:operator_clt}
\cM_{\para} h(\theta):=\sum_{z\in\Lambda}\pi_z(\theta)e^{
\langle \para,z-v\rangle}h(\tau^zF(\theta))=\sum_{z\in\Lambda}e^{\langle \para,z-v\rangle}S_z h.
\end{equation}
Then,
\[
\bE_\nu\left(e^{\frac{i}{\sqrt
N}\langle \para,\Delta_{k}-v\rangle}\;\big|\;\cF_{k}\right)=(\cM_{it/\sqrt{N}} 1)(\tau^{X_{k}}\theta^{k}),
\]
and the reader can then check, by induction, the formula
\begin{equation}\label{eq:clt_formula}
\bE_\nu\left(e^{\frac{i}{\sqrt
N}\langle \para,\tilde X_N\rangle}\right)=\nu(\cM_{it/\sqrt{N}}^N 1).
\end{equation}
The operator $\cM_{\para}'$ acting on the space $\cB$ is an analytic
perturbation of the operator $S'=\cM_0'$. Unfortunately, $S'$ does not
have a nice spectrum on $\cB$, so in order to apply usual perturbation
theory, it is necessary to lift the dynamics to an appropriate space
in the spirit of \cite{BGK}. We do so in section \ref{sec:pert} where we prove the following result.

\begin{lemma}
\label{lem:up}
Under the assumptions of Theorem \ref{thm:mixing} there exists 
$\Ctwo>0$ and a function $\alpha_\para$
analytic near $\{\|\para\|\leq \Ctwo\}$ such that for each  
$n\in\bN$, probability measure $\nu\in\cB$ and local function $f$ depending on $L$ variables we have
\[
\begin{split}
&|\nu({\cM}_{\para}^nf)|\leq CL|\alpha_\para^n|\,|f|_\infty\|\nu\|\\
&\nu(\cM_{\para}^n1)=\alpha_\para^n(1+\cO(\para\|\nu\|))+\cO(\eta^n\|\nu\|).
\end{split}
\]
Moreover, $\alpha_0=1,$  $\dot{\alpha}_0=0$ and $\ddot{\alpha}_0\geq 0$
(the ``dot" stands for the derivatives with respect to $\para$). Finally,  if Assumption \ref{ass:ellipticity} is also satisfied, then
$\ddot{\alpha}_0> 0$.
\end{lemma}
Using Lemma \ref{lem:up} and setting $\Var:= \ddot{\alpha}_0$, we have
\[
\nu(\cM_{i\para/\sqrt N}^N 1)=\alpha_{i\para/\sqrt N}^N(1+\cO(\para N^{-\frac 12}\|\nu\|)+ \cO(\eta^N\|\nu\|).
\]
We can finally compute, for $\|\para\|\leq CN^{\frac 16}$ and $N$ large enough,\footnote{If 
$\|t\|\geq CN^{\frac 16}$ the last statement of Proposition \ref{lem:clt-averaged} is obvious: the left hand side is $\leq 2$.}
\[
\begin{split}
\bE_\nu\left(e^{\frac{i}{\sqrt
N}\langle \para,\tilde X_N\rangle}\right)&=\alpha_{i\para/\sqrt N}^N+\cO\left(\frac{1+\|\para\|}{\sqrt N}\|\nu\|\right)\\
&=\left(1-\frac 1{2N}\langle\para,\Var\para\rangle+\cO(\|\para\|^3N^{-\frac 32})\right)^N+\cO\left(\frac {1+\|\para\|}{\sqrt N}\|\nu\|\right)\\
&=e^{-\frac 1{2}\langle\para,\Var\para\rangle+\cO(\|\para\|^3N^{-\frac 12})}+\cO\left(\frac {1+\|\para\|}{\sqrt N}\|\nu\|\right),
\end{split}
\]
from which Proposition \ref{lem:clt-averaged} follows.
\end{proof}
Next, we need a large deviations estimate.
\begin{lemma}\label{lem:largedev}
Under the assumptions of Theorem \ref{thm:mixing} and Assumption \ref{ass:ellipticity}
there exists $a_0>0$ such that for each $\nu\in \cB$, $n, m\in\bN$ 
and $a\in(0,a_0)$ 
the following holds true
\[
{\bf P}_\nu\left(\left\{\left|\frac 1m(\tilde X_{n+m}-\tilde X_n)\right|\geq
a\right\}\right)\leq Ce^{-Ca^2m}(\|\nu\|+1).
\]
\end{lemma}
\begin{proof}
Again this large deviation result can be obtained by perturbation
theory of the operator $S'$. Indeed, for each $w\in\bR^d$, $\|w\|=1$
and $t\in\bR$, 
\[
\begin{split}
{\bf P}_\nu&\left(\left\{\frac 1m\langle w,\tilde X_{n+m}-\tilde X_n\rangle\geq
a\right\}\right)\leq \bE_\nu\left(e^{\para(\langle w,\tilde X_{n+m}-\tilde X_n\rangle-am)}\right)\\
&\quad\quad=e^{-tam}[(S')^n\nu](\cM_{tw}^m1). 
\end{split}
\]
Since by Theorem \ref{thm:mixing} $\sup_{n\in\bN}\|(S')^n\nu\|\leq C\|\nu\|$, we can apply Lemma \ref{lem:up} and obtain, for $t\leq \Ctwo$,
\[
\begin{split}
{\bf P}_\nu&\left(\left\{\frac 1m\langle w,\tilde X_{n+m}-\tilde X_n\rangle\geq
a\right\}\right)\leq Ce^{-tam}\alpha_{tw}^m(1+C|t|\|\nu\|)+C\eta^m\|\nu\|\\
&\leq Ce^{-tam}\left(1+\frac {t^2}2\langle w,\Var w\rangle+\cO(|t|^3)\right)^m(1+\|\nu\|)+C\eta^m\|\nu\|\\
&\leq Ce^{-tam+\frac {mt^2}2\langle w,\Var w\rangle+\cO(m|t|^3)}(1+\|\nu\|)+C\eta^m\|\nu\|. 
\end{split}
\]
Finally, choosing $\para =\frac {a}{\langle w,\Var w\rangle}$ and $a_0$ so small that the term $\cO(t^3)$ is small with respect to $\frac {a_0^2}{2\langle w,\Var w\rangle}$,  $\frac {a_0}{\langle w,\Var w\rangle}< \Ctwo$ and $e^{-\frac {a_0^2}{2\langle w,\Var w\rangle}}\geq \eta$
\[
{\bf P}_\nu\left(\left\{\frac 1m\langle w,\tilde X_{n+m}-\tilde X_n\rangle\geq
a\right\}\right)\leq Ce^{-Ca^2m}(1+\|\nu\|).
\]
We then conclude by noticing that the above estimate for all the $w$ in the set $\{\pm e_i\}_{i=1}^d$, where $\{e_i\}_{i=1}^d$ is the standard base of $\bR^d$, implies the Lemma.
\end{proof}

\subsection{Main result: Quenched C.L.T}\label{subsec:general}

Let ${\bf P}_\theta$ be the measure ${\bf P}^e$ conditioned to starting the environment in the configuration $\theta$. We will use $\bE_\theta$ for the expectation with respect to ${\bf P}_\theta$.

\begin{thm}
\label{thm:main}
For each dynamics $F$ satisfying assumption \ref{ass:assone} and transition probabilities satisfying assumptions \ref{ass:two} and \ref{ass:ellipticity}, if assumption \ref{ass:pert} is satisfied for $\ve_1$ (where $\ve_0\geq \ve_1>0$ is as in Lemma \ref{lem:variation}), then (using the same notations as in Proposition \ref{lem:clt-averaged}), for $\mue$ almost all $\theta\in\Theta$ the following holds
\begin{enumerate}[(a)]
\item $\frac {1}N X_N \to v$ \;\;\;$ {\bf P}_\theta\; a.s.;$
\item[]
\item $\frac{X_N-vN}{\sqrt{N}}\Rightarrow 
\cN\left(0, \Var\right)$  under ${\bf P}_\theta.$
\end{enumerate}
\end{thm}
\begin{proof}
Lemma \ref{lem:largedev} implies the bound ${\bf P}^e(\{|N^{-1} X_N-v|\geq \ve\})\leq Ce^{-C\ve^{2}N}$; (a) follows then by applying Borel-Cantelli.

To prove (b) let $\alpha\in (0,1)$ be a number to be specified later. 
Combining Lemma \ref{lem:largedev} and Borel-Cantelli
Lemma we see that for any $\delta>0,$ ${\bf P}^e$-almost surely for any $k$ and $0\leq j\leq 2^{(1-\alpha) k}$ 
we have\footnote{The $X$ in $C_{\theta, X, \delta}$ stands for the dependence on the random walk $\{X_n\}_{n\in\bN}$.}
\begin{equation}
\label{eq:smallosc}
\max_{m\in [2^k+j 2^{\alpha k}, 2^k+(j+1) 2^{\alpha k}]} \left|\tX_m-\tX_{2^k+j 2^{\alpha k}}\right|\leq
C_{\theta, X, \delta}\,  2^{(\frac{\alpha}{2}+\delta)k}
\end{equation}
By Fubini Theorem for almost every $\theta$ \eqref{eq:smallosc} holds ${\bf P}_\theta$ almost surely.
Therefore it is enough to prove the convergence along the subsequence $n_{jk}=2^k +j 2^{\alpha k}$
where $0\leq j\leq 2^{(1-\alpha) k}.$

 To conclude it suffices to prove that there exists $\beta>0$ and $b\in\bN$ such that for each smooth 
 function $\vf:\bR^d\to\bR$  compactly supported in a box of size $L$ the following inequality holds 
\begin{equation}\label{eq:hot0}
\bE\left(\left|\bE_\theta(\vf(N^{-\frac 12}\tilde X_N))-\bE_{\cN(0,\Var)}(\vf)\right|^2\right)\leq C_L |\vf|_{\cC^b}N^{-\beta},
\end{equation}
where $\bE_{\cN(0,\Var)}$ is the expectation with respect to the Gaussian measure $\cN(0,\Var)$.
Indeed, denote
$$\xi_{jk}=\bE_\theta\left(\vf\left(\frac{\tX_{n_{jk}}}{\sqrt{n_{jk}}}\right)\right).$$
Then \eqref{eq:hot0} and Chebyshev inequality imply
\begin{equation}\label{eq:hot2}
{\bf P}^e\left(\left\{\left|\xi_{jk}-\bE_{\cN(0,\Var)}(\vf)\right|\geq \ve\right\}\right)\leq C_L|\vf|_{\cC^{b}}\ve^{-2}n_{jk}^{-\beta}.
\end{equation}
Hence, by finally choosing $\alpha$ such that $\alpha+\beta>1$,  $\xi_{jk}\to\bE_{\cN(0,\Var)}(\vf)$ almost surely. 
Next, choose a family $\vf_m$ which is dense in $\cC_0^0(\bR^d).$ Then, for almost every $\theta$, we have
$$ \bE_\theta\left(\vf_m\left(\frac{\tX_{n_{jk}}}{\sqrt{n_{jk}}}\right)\right)\to \bE_{\cN(0,\Var)}(\vf_m) $$
for all $m.$ Then, for any such $\theta$
$$ \bE_\theta\left(\vf\left(\frac{\tX_{n_{jk}}}{\sqrt{n_{jk}}}\right)\right)\to \bE_{\cN(0,\Var)}(\vf) $$ 
for any continuous compactly supported $\vf$ proving (b).

The result is then proved provided \eqref{eq:hot0} is true. It turns out that \eqref{eq:hot0}
can be conveniently interpreted in terms of two independent walks $X_N,Y_N$ in the same environment. In fact, calling $\bE^2$ the expectation with respect to such a process it follows
\[
\begin{split}
\bE&\left(\left|\bE_\theta(\vf(N^{-\frac 12}\tilde X_N))-\bE_{\cN(0,\Var)}(\vf)\right|^2\right)=
\bE^2(\vf(N^{-\frac 12}\tilde X_N)\vf(N^{-\frac 12}\tilde Y_N))\\
&\quad-2\bE(\vf(N^{-\frac 12}\tilde X_N))\bE_{\cN(0,\Var)}(\vf)+\bE_{\cN(0,\Var)}(\vf)^2\\
=&\bE^2(\vf(N^{-\frac 12}\tilde X_N)\vf(N^{-\frac 12}\tilde Y_N))-\bE_{\cN(0,\Var)}(\vf)^2+\cO(L^d\|\vf\|_{\cC^{d+1}}N^{-\beta}).
\end{split}
\]
where we have used the quantitative estimate in the Proposition \ref{lem:clt-averaged}.\footnote{If $\hat\vf$ is the Fourier transform of $\vf$, then Proposition \ref{lem:clt-averaged} yields, for each $\rho<\frac16$ and $b\in\bN$,
\[
\begin{split}
&\left|\bE(\vf(N^{-\frac 12}\tilde X_N))-\bE_{\cN(0,\Var)}(\vf)\right|\leq
C\int_{\|t\|\geq N^{\rho}} |\hat \vf (t)|dt+C\int_{\|t\|\leq N^{\rho}}(1+\|t\|^3)N^{-\frac 12}|\hat\vf(t)|dt\\
&\leq C_dL^d\|\vf\|_{\cC^b}\int_{N^\rho}^\infty x^{-b+d-1}dx+C_dL^d\|\vf\|_{\cC^0}N^{-\frac 12}\int_0^{N^{\rho}} (1+x^3)x^{d-1}dx,
\end{split}
\]
which gives the advertised result provided $b>d$ and $\rho$ is chosen small enough.
}

We have thus reduced the proof of the theorem to proving the following.
\begin{lemma}
\label{lm:hot}
\begin{equation*}
\left|\bE^2(\vf(N^{-\frac 12}\tilde X_N)\vf(N^{-\frac 12}\tilde Y_N))-\bE_{\cN(0,\Var)}(\vf)^2\right|\leq C\|\vf\|_{\cC^b}N^{-\beta}.
\end{equation*}
\end{lemma}
Lemma \ref{lm:hot} is proved in Section \ref{sec:clt-proof}.
\end{proof}

\section{Proofs: Environment}\label{sec:proof_env}
In this section we establish all the needed properties of the
environment dynamics. The basic idea is to prove that the environment enjoys very strong mixing properties.
Our first aim is to prove Theorem \ref{thm:mixing}, that is exponential decay of space-time correlations.

Note that $S1=1$, hence 
\begin{equation}\label{eq:contra}
|S'\mu|\leq |\mu|,
\end{equation}
that is $S$ is a contraction in the $|\cdot|$ norm. In addition, it is possible to prove (although we will not use it here) that there exists $A,B>0$ and $\sigma\in (2\lambda ^{-1},1)$ such that, for each $n\in\bN$ and $\mu\in\cB$,\footnote{See \eqref{eq:banach} for the definition of $\cB$.}
\[
\|(S')^n\mu\|\leq A\sigma^{n}\|\mu\|+B|\mu|.
\]
For finitely many sites the above estimate would suffice to prove that the operator $S'$ is quasi-compact and this, together with the topologically mixing assumption, would imply the existence of a spectral gap. Unfortunately, such a proof is based on the compactness of the unit ball $\{\mu\in \cB\;:\;\|\mu\|\leq 1\}$ in the topology of the $|\cdot|$ norm which fails when one considers infinitely many sites.

The obvious idea is to use explicitly the fact that the dynamics in different sites are independent, hence the system has a product structure, yet this is a subtle issue.
To understand better the situation, let  us recall few fact about the single site systems. At each site we have the dynamical system $(I=[0,1], T)$. Let us consider the norm in $\cM(I)$ given by
\[
\|\nu\|_0=\sup_{|\vf|_{\cC^0}\leq 1}\mu(\vf'),
\]
where $\vf'$ is the derivative of $\vf$.
The Banach space $B=\{\nu\in\cM\;:\;\|\nu\|_0<\infty\}$ consists of measures absolutely continuous with respect to the Lebesgue measure $m_{\cL}$. In addition, if $d\nu=h\, dm_{\cL}$, then the density $h$ is a function of bounded variation and $|h|_{BV}=\|\nu\|$.\footnote{The equality of the norms follows from the usual weak definition of $BV$, the fact that the measures must be absolutely continuos can be easily proved by approximating a measure with finite norm by one with a smooth density (just use a mollifier) and remembering that the 
unit ball of $BV$ is compact in $L^1$. See \cite{KL0} for more details.}
By a change of variable one can compute that, 
if $d\nu=h\, dm_{\cL}$, then $d(T'\nu)=(\cL h) dm_{\cL}$, 
where the operator $\cL$ is defined as
\[
\cL h(x)=\sum_{y\in T^{-1}(x)}|D_yT|^{-1}h(y).
\]
The operator $\cL$ is often called the Ruelle-Perron-Frobenious transfer operator.
It is well known that, if $T$ is topologically mixing, then the operator $\cL$, acting on $BV$ has $1$ as a simple eigenvalue (corresponding to the unique invariant measure absolutely continuous with respect to Lebesgue) and enjoys a spectral gap, that is there exists $\eta_0\in (0,1)$ such that the rest of the spectrum is strictly contained in a disk of radius $\eta_0$ (see, e.g., \cite{Babook} for details).
Clearly the above implies that $T'$ has a spectral gap when acting on $B$. Unfortunately, it turns out that the tensor products of $B$ is a too small a space to be really useful for our purposes.\footnote{The problem is already present for two sites since $BV(I)\otimes BV(I)\neq BV(I^2)$.} This is the reason why we have introduced the spaces $\cB$ which is a generalization of measures with density of bounded variations to the infinite dimensional setting. Yet on such a space $S'$ does not behave very well and we will use an abstract covering space on which the dynamics will exhibit a spectral gap.

More precisely, we would like to introduce a Banach space $\ocB$  and two (possibly only partially defined) maps $\Psi:\cB\to\ocB$ and $\Pro:\ocB\to\cB$ and an operator $\bf S:\ocB\to\ocB$ such that the dynamics of the latter covers the dynamics of $S'$ as illustrated by the following commutative diagram
\begin{equation}
  \label{eq:commutation}
\begin{CD}
    \ocB@> {\bf S}^n >>\ocB\\
     @A{\Psi}AA                         @VV{\Pro}V\\
     \cB @>>{(S')^n}> \cB
  \end{CD}
\end{equation}
We will first define the space $\ocB$ and the map $\bf S$. Then we will prove Theorem \ref{thm:mixing} by proving that $\bf S$ has a spectral gap on $\ocB$. Next we will obtain others, more refined, results by using the same strategy (albeit applied to different operators).

\subsection{Covering dynamics.}\label{sec:covering}
First we define the above mentioned abstract space. Let
$\ocB:=\bC\times\left[\times_{p\in\bZ^d}\cB_p\right]$
where\footnote{For example, if $\nu_p,\nu'_p$
$(p\in\bZ^d)$ are probability measures on $I$ such that $\nu_p=\nu_{p}'$ for all $p\neq
  q$, and we set $\nu:=\otimes_{p\in\bZ^d}\nu_p$,
  $\nu':=\otimes_{p\in\bZ^d}\nu_p'$, then $\nu-\nu'\in\cB_q$.}
\[
\cB_p:=\{\mu\in\cB :\mu(\vf)=0\; \forall
\,\vf\in\cC^0(\Theta)\text{ that do not depend on }\theta_p\}. 
\]
The vector space $\ocB$ is a Banach space when equipped with the norm
\begin{displaymath}
  \|\bmu\|:=\sup\{|c_\mu|,\|\mu_p\|\;:\;p\in\bZ^d\}\ .
\end{displaymath}
Here we use the notational convention that an element $\bmu\in\ocB$ has
components $c_\mu\in\bC$ and $\bar\mu:=(\mu_p)_p$ with
$\mu_p\in\cB_p$.

Next we define a projection $\Pro:D\subset\ocB\to\cB$ and a map
$\Psi:\cB\to\ocB$ allowing to transfer objects between the two spaces.

Let $\mu_*\in \cB$ be a fixed probability measure on $I$ 
then $m=\otimes_{\bZ^d}\mu_*$ is a product probability measure on
$\Theta$. For each $\bmu=(c_\mu,(\mu_p)_p)\in\ocB$ and local
function $f$ we define
\begin{equation}
  \label{eq:decomp}
  \Pro{\bmu}(f):=c_\mu\,m(f)+\sum_{p\in\bZ^d}\mu_p(f),
\end{equation}
which makes clear in which sense $\ocB$ ``covers" $\cB$ ( or $\cM(\Theta)$).

\begin{rem}
  Note that, although $\Pro{\bmu}(f)$ is well defined on each local
  function, $\Pro{\bmu}$ is not necessarily a measure.
  Let $\ocB_M\subset \ocB$ be such that the elements of $\Pro \ocB_M$ give rise to bounded linear functionals on the space of local functions, and hence identify uniquely a measure.\footnote{Since the local functions are dense in the continuous ones by the Stone-Weierstrass theorem.}
 We will call such a measure $\Pro{\bmu}$.
\end{rem}

The choice of the map $\Psi$ is quite arbitrary, we will fix a convenient one.
Consider a strict total ordering $\prec$ of $\bZ^{d}$ such that $0\prec p$ for
each $p\in\bZ^{d}\setminus\{0\}$  and the sets $\{q\;:\;q\prec
p\}$ are finite for each $p\in\bZ^d$.\footnote{For example, one can
start from zero and spiral out over larger and larger cubical shells.} 

Let $q_+$ be the successor of $q$ (that
is, $q\prec q_+$ and there are no $q'\in\bZ^d$ such that $q\prec q'\prec q_+$).
For each $q\in\bZ^d$ we can then consider the $\sigma$-algebra $\cF_q^0$
determined by all the variables $\theta_{q'}$ with $q\preceq q'$, hence $\cF_0^0$
is the complete $\sigma$-algebra. Next, for each $f\in\cC^0(\Theta)$ and
$q\in\bZ^d$, define the operator
$J_qf=m(f\;|\;\cF_q^0)-m(f\;|\;\cF_{q_+}^0)$.  For each local function $f$ we can write\footnote{As $f$ is local
there exists a box $\Lambda_*\subset \Theta$ such that $f$ depends only
on the variables $\{\theta_q\;:\;q\in\Lambda_*\}$, but this means that
the sum consists only of finitely many terms.}
\[
f=m(f)+\sum_{q\in\bZ^d}J_q(f)
\]
Accordingly, for each $\mu\in\cB$ we define
$\mu_q(f):=J_q'\mu(f)\in\cB_q$, and the
lift
\begin{displaymath}
  \Psi(\mu):=(\mu(1),(J_q'\mu)_q)\ .  
\end{displaymath}
Note that $\Psi$ is a bounded operator. Indeed if $q\prec p$, then
\[
|J_q'\mu (\partial_{\theta_p}\vf)|\leq |\mu|\, |J_q(\partial_{\theta_p}\vf)|_\infty\leq 2 |\mu|\,\|\mu_*\|\,|\vf|_\infty.
\]
If $q\succ p$, then
\[
|J_q'\mu (\partial_{\theta_p}\vf)|\leq |\mu(\partial_{\theta_p}J_q\vf)|\leq \|\mu\|\, |J_q(\vf)|_\infty\leq 2 \|\mu\|\,|\vf|_\infty.
\]
Finally, for $q=p$, we have $|J_q'\mu (\partial_{\theta_p}\vf)|\leq (|\mu|\,\|\mu_*\|+\|\mu\|)|\vf|_\infty$. In other words there exists $\Cseven>0$, depending on the choice of $\mu_*$, such that
\begin{equation}\label{eq:J-bound}
\|\Psi(\mu)\|\leq \Cseven \|\mu\|.
\end{equation}
Clearly, $\Psi(\cB)\subset D$ and for each $\mu\in\cB$ it holds true
\[
\Pro(\Psi(\mu))=\mu.
\]
Now that we know how to lift measures, we can address the
dynamics. 

For all $z\in\bZ^d$, $\tau^z\mu_q(\vf):=\mu_q(\vf\circ
\tau^z)=0$ if $\vf$ does not depend on $\theta_{q-z}$. Thus we can define
the decomposition for $\tau^{z}\mu$ via the decomposition $\mu_q=J_q'\mu$ of $\mu$:
\[
\tau^{z}\mu=\sum_{q\in\bZ^d}(\tau^{z}\mu)_q:=\sum_{q\in\bZ^d}\tau^{z}\mu_{q-z}.
\]
Setting $\Lambda_1=\cup_{z\in\Lambda}\tau^{-z}\Lambda$,  we can define the covering dynamics
${\bf S}{\bmu}$ by $(c_\mu, \bar S\bar\mu+\bar \zeta c_\mu)$ where
$\bar\zeta:=(\zeta_q)$ with $\zeta_q:=J_q'S'm=:J_q'\zeta$ and 
\begin{equation}\label{eq:dynone}
(\bar S\bar\mu)_p=\begin{cases}\sum_{z\in\Lambda}S_z'\mu_{p-z}+\sum_{q\in\Lambda_1}\sum_{z\in\Lambda}J_p'\hat S_z'\mu_{q-z}&\text{ for each }p\not\in \Lambda_1,\\
\sum_{z\in\Lambda}A_z'\mu_{p-z}+\sum_{q\in\Lambda_1}\sum_{z\in\Lambda}J_p'\hat S_z'\mu_{q-z}&\text{ for each }p\in \Lambda_1,
\end{cases}
\end{equation}
where $S_zf(\omega):=\pi_z(\omega)f\circ F\circ \tau^z(\omega)$, $A_zf(\omega):=a_z f\circ F\circ \tau^z(\omega)$ and $\hat S_z:=S_z-A_z$.
It is easy to check the following.
\begin{lemma}\label{lem:covering}
The operator $\bf S$ is well defined as a bounded operator from $\ocB$ to $\ocB$.
For each $\bmu\in\ocB$ and continuous local function $f$ we have $\Pr({\bf S}\bmu)(f)=\Pr(\bmu)(Sf)$ which
implies that for each $n\in\bN$ and 
$\mu\in\cB$ we have  $\Pr {\bf S}^n\Psi\mu={S^n}'\mu$
\end{lemma}

We have thus established a setting in which the commutative diagram \eqref{eq:commutation} holds true.

\subsection{Mixing properties of the environment}\label{sec:mixone}
We have now the necessary machinery to deal with the statistical properties of the environment.

\begin{proof}[{\bfseries Proof of Theorem \ref{thm:mixing} }]
Let us first discuss the environment dynamics $F$. Its basic properties are described by the so called Lasota-Yorke inequalities asserting that there exists $B>0$ such that, for all $n\in\bN$,\footnote{The first is trivial since $|F'\mu(\vf)|=|\mu(\vf\circ F)|\leq|\mu|\,|\vf|_\infty$.
For the second, given any smooth local function $\vf$, let $\vf_{\theta_{\neq q}}(\xi):=\vf(\theta^{\xi})$, where $\theta^{\xi}_p=T\theta_p$ for each $p\neq q$ while $\theta^{\xi}_q=\xi$. Next introduce a function $\phi$, piecewise linear in the variable $\theta_q$ such that $\vf_{\theta_{\neq q}}(\theta_p)-\phi(\theta)=0$ for each $\theta_p$ on the discontinuity values of $T$. By construction $\vf_{\theta_{\neq q}}(T\theta_p)-\phi(F\theta)$ is then a Lipschitz function in $\theta_p$, thus
\[
\begin{split}
|F'\mu(\partial_{\theta_q}\vf|&=|\mu((\partial_{\theta_q}(\vf-\phi))\circ F)|+
|\mu|\,|\partial_{\theta_q}\phi|_\infty
\leq |\mu(\partial_{\theta_q}(|D_{\theta_q}T|^{-1}(\vf-\phi)\circ F)|+
C|\mu|\,|\vf|_\infty\\
&\leq \|\mu\|\, |D_{\theta_q}T|^{-1}(\vf-\phi)\circ F|_\infty+C|\mu|\,|\vf|_\infty
\leq \left[2\lambda^{-1} \|\mu\|+C|\mu|\right]|\vf|_\infty.
\end{split}
\]
The above yields $\|F'\mu\|\leq 2\lambda^{-1} \|\mu\|+C|\mu|$ which iterated yields the wanted result with $B=(1-2\lambda^{-1})^{-1} C$. See \cite{Babook, KL0} if more details are needed.}
\begin{equation}\label{eq:lasota_0}
\begin{split}
&|F'\mu|\leq |\mu|\\
&\|(F')^n\mu\|\leq (2\lambda^{-1})^n\|\mu\|+B|\mu|.
\end{split}
\end{equation}
Note that the above implies that $\{\|(F')^n\mu\|\}_{n\in\bN}$ is bounded. 
\begin{lemma}\label{lem:gap0}
There exists
$\eta_*\in (\eta_0,1)$ such that, for each $q\in\bZ^d$, $\mu_q\in\cB_q$,
\[
\|(F')^n\mu_q\|\leq C\eta_*^n\|\mu_q\|.
\]
\end{lemma}
\begin{proof}
For each local function $\vf\in\cC^0$, we can define
$\vf_{\theta_{\neq q}}^n(\xi):=\vf(\theta^{n,\xi})$, where $\theta^{n,\xi}_p=T^n\theta_p$ for each $p\neq q$ while $\theta^{n,\xi}_q=\xi$,
\[  
\begin{split}
|(F')^n\mu_q(\vf)|&=\left|\mu_q\left(\partial_{\theta_q}\int_0^1[\chi_{[0,\theta_q]}(\xi)-\theta_q] \vf_{\theta_{\neq q}}^n(T^n\xi) d\xi\right)\right|\\
&\leq \|\mu_q\|\cdot\left|\int_0^1\cL^n[\chi_{[0,\theta_q]}(\xi)-\theta_q] \cdot \vf_{\theta_{\neq q}}^n(\xi) d\xi\right|_\infty\\
&\leq C\|\mu_q\|\,\left|\cL^n[\chi_{[0,\theta_q]}-\theta_q]\right|_{BV}|\vf|_\infty\leq C\|\mu_q\|\eta_0^n|\vf|_\infty,
\end{split}
\]
by the spectral gap of $\cL$ and the fact that $m_{\cL}(\chi_{[0,\theta_q]}-\theta_q)=0$, i.e. it is a zero average function.
Then, by the Lasota-Yorke inequality,
\[
\begin{split}
\|(F')^{j+k}\mu_q\|&\leq (2\lambda^{-1})^j\|(F')^{k}\mu_q\|+B|(F')^{k}\mu_q|\\
&\leq (2\lambda^{-1})^j\left[(2\lambda^{-1})^k\|\mu_q\|+B|\mu|\right]+B C\eta_0^k\|\mu_q\|\\
&\leq \left[(2\lambda^{-1})^{j+k}+ (2\lambda^{-1})^jB+BC\eta_0^k\right]\|\mu_q\|.
\end{split}
\]
The result follows by optimizing the choice of $j+k=n$.
\end{proof}

\begin{lemma}
\label{SLMult}
Multiplication by a $\cC^1$ local function is a bounded operator on $\cB.$ 
\end{lemma}
\begin{proof}
For any smooth local functions $\psi,\phi$, $|\psi|_\infty\leq 1$ we have
\begin{equation}\label{eq:times}
|\nu(\phi\cdot \partial_{\theta_i}\psi)|=\left|\nu\left(\partial_{\theta_i}\int_0^{\theta_i}(\phi\partial_{\theta_i}\psi)\right) \right | \leq \|\nu\|\left|\int_0^{\theta_i}\phi\partial_{\theta_i}\psi\right|_\infty\leq 3\|\nu\|\, |\phi|_{\cC^1}.
\end{equation}
\end{proof}
To use the above facts, it is convenient to introduce a more compact notation for the pieces that make up the operator $\bar S$. Let  $\Id_A:\bZ^d\to \{0,1\}$ be the characteristic function of the set $A\subset\bZ^d$. Then define the operators $K_{z,p,q,\sigma}:\cB\to\cB$ by $K_{z,p,q,0}:=\Id_{\{p\}}(q+z)A_z'$ and $K_{z,p,q,1}:=\Id_{\Lambda_1^c}(p)\Id_{\{p\}}(q+z) \hat S'_z+\Id_{\Lambda_1}(q+z)J_p'\hat S'_z$. With this notation \eqref{eq:dynone} can be rewritten as
\[
(\bar S\bar\mu)_p=\sum_{z\in\Lambda}\;\sum_{\sigma\in\{0,1\}}\;\sum_{q\in\bZ^d}K_{z,p,q,\sigma}\mu_{q}.
\]
Hence, iterating,
\begin{equation}\label{eq:Smu}
(\bar S^n\bar\mu)_{q_0}=\!\!\!\!\sum_{z_1,\dots, z_n\in\Lambda}\;\sum_{\sigma_1,\dots,\sigma_n\in\{0,1\}}\;\sum_{q_1,\dots, q_n\in\bZ^d}\!\!\!\! K_{z_1,q_0,q_1,\sigma_1}\cdots K_{z_n,q_{n-1},q_n,\sigma_n}\mu_{q_n}.
\end{equation}
By Assumption \ref{ass:pert}, Lemma \ref{lem:gap0}, Lemma \ref{SLMult} and the inequalities \eqref{eq:lasota_0} it follows that there exists a constant $\Csix>0$, depending only on $F$ and $\pi_z$, such that
$\sum_q\|K_{z,p,q,1}\mu_q\|\leq \Csix\, \ve a_z\|\bar \mu\|$ and 
\[
\sum_{q_1,\dots, q_n\in\bZ^d}\|K_{z_1,q_0,q_1,0}\cdots K_{z_n,q_{n-1},q_n,0}\mu_{q_n}\|\leq \Csix \, a_{z_1}\cdots a_{z_n}\eta_*^n\|\bar\mu\|
\]
Accordingly, if $\Csix^2\ve+\eta_*<1$, then there exists $n_*\in\bN$ and $\eta\in(\eta_*,1)$ such that $\Csix(\eta_*+\Csix^2\ve)^{n_*}\leq \eta^{n_*}<1$. This means that every $z_1,\dots,z_{n_*}$ term in \eqref{eq:Smu} will be smaller than $\eta^{n_*}a_{z_1}\cdots a_{z_{n_*}}\|\bar\mu\|$, hence
for all $n\in\bN$, 
\begin{equation}\label{eq:gap}
\| (\bar S^{n}\bar\mu)\|\leq\Csix\,  \eta^{n-n_*}\|\bar\mu\| \,.
\end{equation}
Since ${\bf S}^n\bmu=(c_\mu, \bar S^n\bar \mu+c_\mu\sum_{k=1}^{n-1}\bar S^k\bar\zeta )$ and the series $\bar\zeta_*=\sum_{k=1}^\infty \bar S^k\bar\zeta$ converges by \eqref{eq:gap}, it follows that $\bmuinv:=(1,\bar\zeta_*)$ is an invariant vector for ${\bf S}$. In addition
\begin{equation}\label{eq:gap2}
\|{\bf S}^n\bmu-c_\mu\bmuinv\|\leq C\eta^{n}\|\bmu\|,
\end{equation}
That is the operator $\bf S$ on $\ocB$ has one as a simple maximal eigenvalue and a spectral gap.
From this result we can obtain the decay of temporal correlation simply by projecting down to $\cB$. Indeed, let $\mu$ be a probability measure and $\phi$ be a smooth local function depending only on the sites $A\subset \bZ^d$ and let $L$ be the cardinality of $A$, then, by Lemma \ref{lem:covering} and \eqref{eq:gap2}, \eqref{eq:J-bound}, 
\begin{equation}\label{eq:equi}
\begin{split}
|\mu(\phi\circ S^n)&-\Pro(\bmuinv)(\phi)|= |\Pro({\bf S}^n(\Psi(\mu)-\bmuinv))(\phi)|\\
&=\left|\sum_{q\in A} ({\bf S}^n(\Psi(\mu)-\bmuinv))_q(\phi)\right|\\
&\leq\sum_{q\in A} \|({\bf S}^n(\Psi(\mu)-\bmuinv))_q\| \,|\phi|_\infty\leq CL\eta^n(\|\mu\|+C)|\phi|_\infty.
\end{split}
\end{equation}
Thus, remembering \eqref{eq:contra}, $\bmuinv\in\ocB_M$, that is it gives rise to a bounded linear functional on local functions. Accordingly, we can define the measure $\muinv=\Pro(\bmuinv)$ which will be invariant by $S'$.
Equation \eqref{eq:equi} gives then the temporal correlation decay for such a measure.

To have the spatial decay of correlations note that if $\varphi$ and $\phi$ are supported at a
distance $M$, then their support, under the dynamics, grows at most
linearly in time, thus it will take a time $\frac{M}{2\Czero}$ before the supports have a common variable. 
Accordingly, since $\vf\phi$ depends on $2L$ variables, \eqref{eq:equi} (applied repeatedly to the product measure $m$) implies
\[
\begin{split}
\muinv(\varphi\phi)&=(S')^{M/2\Czero}m(\varphi\phi)+\cO(L\eta^{M/2\Czero}|\varphi\phi|_\infty)\\
&=m(S^{M/2\Czero}\varphi)\,m(S^{M/2\Czero}\phi)
+\cO(L\eta^{M/2\Czero}|\varphi\phi|_\infty)\\
&=\muinv(g)\muinv(\phi)+\cO(L\eta^{M/2\Czero}|\varphi\phi|_\infty). 
\end{split}
\]
\end{proof}

\subsection{Perturbation Theory}
\label{sec:pert}
In this section we prove Lemma \ref{lem:up}.

We deal with operators of the the type $\cM_{\para}f:=\sum_{z\in\Lambda}S_z(e^{\langle\para,z-v\rangle}f)$
where $\para\in\bC^d$.  The problem is to study the spectrum for small $\para$.\footnote{This problem is already well investigated, see in particular \cite{BGK}, here we treat it in detail only because we need some explicit estimates not readily available in the literature.}

\begin{proof}[{\bfseries Proof of Lemma \ref{lem:up}}]

First of all we need to lift the operator to our covering space. The obvious solution is to define
${\bf M}_{\para}{\bmu}$ by 
\[
\left(c_\mu m(\cM_\para 1)+\sum_{p\in\Lambda_1}\sum_{z\in\Lambda}e^{\langle \para, z-v\rangle}\mu_{p-z}(\hat S_{z}1),\  \bar S_{\para}\bar\mu+\bar\zeta c_\mu\right)
\]
where $\bar\zeta=(\zeta_q):=(J_q' \cM_\para'm)$ and
\begin{equation}\label{eq:dytwo}
(\bar S_{\para}\bar\mu)_q=\begin{cases}\sum_{z\in\Lambda}e^{\langle\para, z-v\rangle}S_{z}'\mu_{q-z}+\sum_{\substack{p\in\Lambda_1\\z\in\Lambda}}e^{\langle\para, z-v\rangle}J_q'\hat S_{z}'\mu_{p-z}&\forall\;\;q\not\in \Lambda_1,\\
\sum_{z\in\Lambda}e^{\langle\para, z-v\rangle}A_z'\mu_{q-z}+\sum_{\substack{p\in\Lambda_1\\z\in\Lambda}}e^{\langle\para, z-v\rangle}J_q'\hat S_{z}'\mu_{p-z}&\forall\;\;q\in \Lambda_1.
\end{cases}
\end{equation}
A direct computation shows that, for each smooth local function $\vf$, $\Pr({\bf M}_\para\bmu)(\vf)=\Pr(\bmu)(\cM_\para\vf)$, thus the lift covers the dynamics. In addition, one can easily check that ${\bf M}_0={\bf S}$ and that ${\bf M}_\para$ is analytic in $\para$.\footnote{For the first assertion note that
\[
\sum_{p\in\Lambda_1}\sum_{z\in\Lambda}\mu_{p-z}(\hat S_{z}1)=
\sum_{p\in\bZ^d}\sum_{z\in\Lambda}\mu_{p-z}(\hat S_{z}1)=\sum_{p\in\bZ^d}\sum_{z\in\Lambda}\mu_{p}(S_{z}1)
=\sum_{p\in\bZ^d}\mu_p(S1)=0.
\]
For the latter just write it as power series of $\para$.}
Accordingly, standard perturbation theory implies that there exists $\alpha_\para$, ${\bmu}_\para$, analytic in $\para$, such that ${\bf M}_\para{\bmu}_\para=\alpha_\para {\bmu}_\para$ with $\alpha_0=1$, $\bmu_0=\bmuinv$. 

We will normalize $\bmu_\para$ so that $\bmu_\para=(1,\bar\mu_\para)$. Setting $\mu_\para:=\Pr (\bmu_\para)$, for each fixed local function $f$, $\mu_\para(f)=\Pr\bmu_\para (f)$ is analytic in $\para$ since the sum implicit in the right hand side is just a finite sum.\footnote{Note that $\mu_\para$ is not necessarily a measure and gives rise to an analytic object only when applied to a local function. We will abuse notations by writing $\dot\mu_\para$ to mean the functional on local functions defined by $\Pro(\frac{d}{d\para}\bmu_t)(\vf)$.} However Lemma \ref{lem:up} requires a more quantitative information.

By the arguments of section \ref{sec:mixone} (see \eqref{eq:gap}) it follows that ${\bf M}_0={\bf S}=\Pi+R$ where $\Pi^2=\Pi$, $\Pi R=R\Pi=0$ and $\|R^n\|\leq C\eta^n$, for all $n\in\bN$. Thus by standard perturbation theory (see \cite{Ka}), ${\bf M}_\para=\alpha_\para\Pi_\para+R_\para$ where $|\alpha_0-\alpha_\para|\leq C\|\para\|$, $\|\Pi_t-\Pi\|\leq C\|t\|$, $\|R_\para^n\|\leq C\eta_\para^n$, with $\eta_\para\leq\eta+C\|t\|$. Hence, $\|{\bf M}_\para^n\|\leq C|\alpha_\para|^n+C\eta_\para^n$ and, for each local function function $f$ depending only on $L$ variables
\[
|\cM_\para^n\nu(f)|=|\Pr({\bf M}_\para^n\Psi(\nu))(f)|\leq L|f|_\infty\|{\bf M}_\para^n\Psi(\nu)\|\leq C|\alpha_\para|^n \|\nu\|L|f|_\infty
\]
provided that $|\alpha_\para|\geq \eta_\para$ which holds for all $\|t\|<B$ for some $B>0$. This proves the first inequality of Lemma \ref{lem:up}.

To prove the second note that $\Pi_\para\bnu=\ell_\para(\bnu)\bmu_\para$ with $\ell_0(\bnu)=[\bnu]_0$ hence\footnote{Here we use the notation $[\bnu]_0$ to designate the components $c_\nu$ of the vector $\bnu=(c_\nu,\bar\nu)$.}
\[
\begin{split}
\cM_\para^n\nu(1)&=\left[{\bf M}_\para^n\Psi(\nu)\right]_0=\left[\alpha_\para^n\Pi_\para\Psi(\nu)+\cO(\eta_\para^n\|\nu\|)\right]_0\\
&=\left[\alpha_\para^n\Pi_0\Psi(\nu)\right]_0+\cO((\eta_\para^n+\alpha_\para^nC t)\|\nu\|)
=\alpha_\para^n(1+\cO(\para \|\nu\|))+\cO(\eta_\para^n\|\nu\|).
\end{split}
\]
 
 Finally, to study the derivatives of $\alpha$ 
 we use the relation $\mu_t(\cM_t\vf)=\alpha_t \mu_t(\vf)$ for any local smooth function $\vf$. 
 Differentiating with respect to $\para$ yields
\begin{equation}\label{eq:pert-one}
\begin{split}
&\dot\mu_\para(\cM_\para\vf)+\mu_\para(\dot {\cM}_\para\vf)=\dot\alpha_\para\mu_\para(\vf)+\alpha_\para\dot\mu_\para(\vf)\\
&\ddot\mu_\para(\cM_\para\vf)+2\dot\mu_\para(\dot { \cM}_\para\vf)+\mu_\para(\ddot \cM_\para\vf)= \ddot\alpha_\para\mu_\para(\vf)+2\dot\alpha_\para\dot\mu_\para(\vf)+\alpha_\para\ddot\mu_\para(\vf)\\
&\dot\mu_\para(1)=\ddot\mu_\para(1)=0,
\end{split}
\end{equation}
Since $\dot\cM_\para=\sum_{z\in\Lambda}(z-v)e^{\langle t,z-v\rangle}S_z$ and $\ddot\cM_\para=\sum_{z\in\Lambda}(z-v)\otimes (z-v)e^{\langle t,z-v\rangle}S_z$ the above equations, for $\para=0$ imply (substituting $\vf=1$)
\[
\dot\alpha_0=\muinv(\dot\cM_0 1)=\muinv (g-v)=0,
\]
where $g$ is defined in \eqref{eq:S} and we have used \eqref{eq:average-one}.
Next, substituting in the first of the \eqref{eq:pert-one}, $\vf=\sum_{k=0}^{n-1}S^k\phi$, for some local function $\phi$, we have
\[
\sum_{k=0}^{n-1}\dot\cM_0'\muinv(S^k\phi)=\dot\mu_0((\Id-S)\sum_{k=0}^{n-1}S^k\phi)=\dot\mu_0(\phi)-
\Pro({\bf S}^{n-1}\dot\bmu_0)(\phi).
\]
By \eqref{eq:gap2}, taking the limit for $n$ to infinity, we have
\begin{equation}\label{eq:mu-der}
\dot\mu_0(\phi)=\sum_{k=0}^{\infty}\dot\cM_0'\muinv(S^k\phi).
\end{equation}
Finally, the second of the \eqref{eq:pert-one}, setting $\vf=1$ and $\para=0$, yields\footnote{Remember that $\ddot \alpha_\para=(\partial_{\para_i}\partial_{\para_j}\alpha_\para)$ is a $d\times d$ matrix.}
\begin{equation}\label{eq:sec-der}
\begin{split}
\ddot\alpha_0&=2\sum_{n=0}^\infty\muinv(\dot\cM_0S^n\dot\cM_01) +\muinv(\ddot \cM_0 1)\\
&=2\sum_{n=0}^\infty\bE_{\muinv}\left(\tilde \Delta_n\otimes\tilde \Delta_0\right)+\bE_{\muinv}\left(\tilde \Delta_0\otimes\tilde\Delta_0\right).
\end{split}
\end{equation}
Since $\dot \cM 1$ is a local function the sum is convergent. Hence
\[
\begin{split}
\ddot{\alpha}_0&=\lim_{n\to\infty}\frac 1n\left[2\sum_{k,m=0}^n\bE_{\muinv}\left(\tilde \Delta_{m+k}\otimes\tilde \Delta_k\right)+\bE_{\muinv}\left(\tilde \Delta_k\otimes\tilde\Delta_k\right)\right]\\
&=\lim_{n\to\infty}\frac 1n\bE_{\muinv}\left(\tX_n\otimes \tX_n\right)\geq 0.
\end{split}
\]

Finally, if there exists $w\in\bR^d$ such that $\ddot{\alpha}_0w=0$, it means (from \eqref{eq:sec-der} and Theorem  \ref{thm:mixing}) that there exists a constant $\Ceight>0$ such that, for all $n\in\bN$,\footnote{The latter follows by Theorem \ref{thm:mixing}. Let $\beta=\sum_z\langle w,z-v\rangle^2\pi_z$, and $G_j=\sum_z \langle w, z-v\rangle S_z(S^j\langle w,g-v\rangle)$, then
\[
\bE\left(|\langle w,\tX_n\rangle|^2\right)=\sum_{k=0}^{n-1}\mue(S^k \beta)+\sum_{j=1}^{n-1}\sum_{k=0}^{n-j-1}\mue(S^kG_j)
\]
and $|\mue(S^k \beta)-\muinv(\beta)|\leq C\eta^{k}$, $|\mue(S^kG_j)-\muinv(G_j)|\leq Cj^d\eta^k$.
Moreover, setting $\nu_{k,z}(\vf):= \langle w, z-v\rangle\mue(S^k S_z\vf)$, for $\ve$ in assumption \ref{ass:pert}
such that $(1+3\ve)2\lambda^{-1}<1$ we have (by equation \eqref{eq:lasota_0} and \eqref{eq:times}) $\|\nu_{k,z}\|\leq C$. Hence Theorem \ref{thm:mixing} yields
\[
|\mue(S^kG_j)|\leq |\sum_z|\nu_{k,z}(S^j(g-v))|\leq C\eta^j
\]
Thus we can write
\[
\left|\bE\left(|\langle w,\tX_n\rangle|^2\right)-\bE_\muinv\left(|\langle w,\tX_n\rangle|^2\right)\right|\leq C\sum_{k=0}^n\eta^k+C\sum_{j=1}^{n-1}\left[\sum_{k=0}^j \eta^j+\sum_{k=j+1}^{n-j}j^d\eta^k\right]\leq C.
\]
}
\[
\bE_\muinv\left(|\langle w,\tX_n\rangle|^2\right)\leq \Ceight\;; \quad \bE\left(|\langle w,\tX_n\rangle|^2\right)\leq \Ceight
\]
We can thus extract a subsequence $\{n_j\}$ such that $\langle w,\tX_{n_j}\rangle$ converges weakly almost surely
to a random variable $Z$. Let $\psi=\bE_\muinv(Z\;|\; \cF_0)$ and $\tgw=\langle w,g-v\rangle$, then, for each $\cF_0$ measurable smooth local function $\vf$,
\[
\begin{split}
\bE_\muinv(\vf (\psi-S\psi))&=\lim_{j\to\infty}\bE_\muinv(\vf(\tX_{n_j}-\tX_{n_j+1}))=\lim_{j\to\infty}\muinv(\vf(\tgw-S^{n_j+1}\tgw))\\
&=\muinv(\vf\tgw),
\end{split}
\]
where we have used Theorem \ref{thm:mixing}. Thus $\tgw=\psi-S\psi$, $\mue$-a.s.. This implies that, setting $M_0=0$, and
\[
\begin{split}
M_{n+1}-M_n&=\langle w,\tDelta_{n}\rangle-\bE(\langle w,\tDelta_n\rangle\;|\;\cF_n)+\psi(\omega^{n+1})-S\psi(\omega^n)\\
&=\langle w,\tDelta_{n}\rangle+\psi(\omega^{n+1})-\psi(\omega^n),
\end{split}
\]
$M_n$ is a $\bP_\muinv$ stationary martingale. Moreover,
\[
\langle w,\tX_n\rangle=M_n-\psi(\omega^n)+\psi(\omega^0).
\]
From this it follows that 
\[
\begin{split}
C\geq\bE_\muinv(|M_n|^2)&=\sum_{k=1}^{n-1}\bE_\muinv(|\langle w,\tDelta_{n}\rangle+\psi(\omega^{n+1})-\psi(\omega^n)|^2)\\
&=\sum_{k=1}^{n-1}\bE_\muinv(|\langle w,\tDelta_{n}\rangle+\psi(\omega^{n+1})|^2-|\psi(\omega^n)|^2)\\
&=(n-1)\left[\bE_\muinv(|\langle w,\tDelta_{1}\rangle+\psi(\omega^{1})|^2-|\tgw(\omega^0)+S\psi(\omega^0)|^2)\right].
\end{split}
\]
Thus $\sum_z\pi_z|\langle w,z-v\rangle+\psi\circ F\circ \tau^z|^2=\left|\sum_z\pi_z(\langle w,z-v\rangle+\psi\circ F\circ \tau^z)\right|^2$, that is $\langle w,z-v\rangle+\psi\circ F\circ \tau^z=\tgw+S\psi=\psi$, $\muinv$-a.s.. 

Next, let $\alpha_z=\muinv(\pi_z)$, then, $\sum_z\alpha_z=1$ and $\sum_z\langle w, z-v\rangle\alpha_z=\muinv(\tgw)=0$.
Hence,
\[
\sum_z\alpha_z \psi\circ F\circ \tau^z=\psi\quad \muinv\text{-a.s.} .
\]
Note that the operator $S_\alpha\vf:=\sum_z\alpha_z \vf\circ F\circ \tau^z$ defines a Markov process with invariant measure $\mue$ and satisfies the hypothesis of Theorem \ref{thm:mixing}.  Since $\psi=\sum_{k=0}^{n-1}S^k\tgw+S^{n}\psi$, for each $\phi\in L^2(\muinv)$, 
\[
\lim_{j\to\infty}\muinv(\phi \sum_{k=0}^{n_j-1}S^k\tgw)=\lim_{j\to\infty}\bE_\muinv(\phi(\omega^0)\langle w,X_{n_j}\rangle)=\muinv(\phi \,\psi).
\]
In addition, assumption \ref{ass:pert} implies that setting,\footnote{Indeed, $\muinv(\pi_z)\leq a_z (1+\ve)$, thus
$\pi_z\geq (1-\ve)a_z\geq(1-\ve)(1+\ve)^{-1}\alpha_z$.} for each smooth local function $\phi$, $\nu_{n,\phi}(\vf)=\muinv(\phi S_\alpha^n\vf)$,
\[
|\nu_{n,\phi}(\vf)|\leq \frac{(1+\ve)^n}{(1-\ve)^n}|\phi|_\infty\muinv (S^n|\vf|)= \frac{(1+\ve)^n}{(1-\ve)^n}|\phi|_\infty\muinv(|\vf|).
\]
Thus $\nu_{n,\phi}$ is absolutely continuos with respect to $\muinv$ with density $\rho_{n,\phi}\in L^\infty(\muinv)$. Accordingly,\footnote{Note that \eqref{eq:lasota_0} imply $\|S_\alpha'\mu\|\leq\sum_z\alpha_z\|F'\mu\|\leq(2\lambda)^{-1}\|\mu\|+B|\mu|$ and $|S_\alpha\mu|\leq|\mu|$. Thus, iterating, for each $n\in\bN$, $\|(S_\alpha^n)'\mu\|\leq C\|\mu\|$.}
\[
\begin{split}
\muinv(\phi\psi)&=\muinv (\phi S_\alpha^n \psi)=\lim_{j\to\infty}\sum_{k=0}^{n_j-1}\muinv (\phi S_\alpha^n S^k \tgw)\\
&=\sum_{k=0}^{n_l-1}\left[\mue( S^k \tgw)\muinv(\phi)+\cO(C_\phi\eta^{n}k^d)\right]
+\lim_{j\to\infty}\sum_{k=n_l}^{n_j-1}\cO(C_\phi\eta^k)\\
&=\sum_{k=0}^{n_l-1}\mue( S^k \tgw)\muinv(\phi)+C_\phi\cO(\eta^{n_l}+\eta^nn_l^d).
\end{split}
\]
Taking first the limit for $n\to\infty$ and the one $l\to\infty$  yields $\muinv(\phi\psi)=\muinv(\phi)\mue(\psi)$. That is $\psi$ is $\muinv$ almost surely constant. This implies that $\tgw=0$ and hence $\langle w,z-v\rangle\pi_z=0$, $\muinv$ a.s.. This is equivalent to saying that the vectors in the set
$\{(\langle e_1,z\rangle,\dots, \langle e_d,z\rangle)\}_{z\in\Lambda}\cup\{(1,\dots,1)\}$ 
are linearly dependent over $\bR$, but this implies that they are 
linearly dependent over $\bZ$. In other words we can assume that 
$w\in\bZ^d$. Finally, since $\pi_z$ is smooth, we have 
$\langle w,z\rangle=\langle w,v\rangle$ unless $\pi_z\equiv 0$, 
which contradicts  Assumption \ref{ass:ellipticity}.
\end{proof}

\subsection{Variation bounds for conditional measures}
\label{sec:variation}
In the previous subsection we obtained several results for random walks provided that we start the environment in a measure with ``density" of bounded variation. 
Here we show why such measures constitute a natural class for the 
problem at hand. More precisely we shall show that if we start with a nice measure and condition on a behavior of a walk during an initial time interval we still have a good control on the variation of densities.

For future needs we consider two random walks $(X_t,Y_t)$ evolving in the same environment starting respectively at $a,b\in\bZ^d$ and with the environment at time zero distributed according to the measure $\nu\in\cB$. Let $\bP_{a,b,\nu}^2$ be the measure on $(\Theta\times \bZ^{2d})^\bN$  associated to such a process and $\bE_{a,b,\nu}^2$ the corresponding expectation $\bP^2_\nu:=\bP^2_{0,0,\nu}$ and $\bE^2_\nu:=\bE^2_{0,0,\nu}$.

Let  $m\in\bN$ and consider the $\sigma$-algebra
$\cF_m^{XY}=\sigma\{X_1,Y_1,\dots,X_m,Y_m\}$. We are interested in
computing $\bE_{a,b,\nu}^2(f(X,Y,\theta^m)\;|\; \cF_m^{XY})$ for each local $\cF_m^{XY}\otimes\cT$-measurable function $f$ and probability measure $\nu\in\cB$. Thus, we are interested in the measures $\nu^{XY}_{a,b,m}$ defined by
\[
\bE^2_{a,b,\nu}(f(X,Y,\theta^m)\;|\; \cF_m^{XY})=:\int_\Theta f(X,Y,\theta)\;\nu^{XY}_{a,b,m}(d\theta).
\]
\begin{lemma}
\label{lem:variation}
There exists $\Cfive>0$ and $0<\ve_1\leq\ve_0$ such that, if assumption \ref{ass:pert} is satisfied for $\ve_1$, then for each $m\in\bN$, $a,b\in \bZ^d$ and probability measure $\nu\in\cB$
the following holds
\[
\|\nu^{XY}_{a,b,m}\|\leq \Cfive\|\nu\|.
\]
\end{lemma}
\begin{proof}
Given two random walks realizations $X,Y:\bN\to\bR^d$, let us define the operators
\[
S_{ X,Y,k}f(\theta):=\pi_{z_k}(\tau^{X_k}\theta)\pi_{w_k}(\tau^{Y_k}\theta)
f\circ F(\theta),
\]
where $z_k=X_{k+1}-X_k$ and $w_k=Y_{k+1}-Y_k$.
With such a notation we can write
\[
\nu^{XY}_{a,b,m}=\frac{S_{X,Y, m-1}'\cdots  S_{X,Y, 0}'\nu}{S_{X,Y, m-1}'\cdots  S_{X,Y, 0}'\nu(1)}.
\]
Recalling \eqref{eq:times} and the Lasota-Yorke inequality for the map $F$ (see \eqref{eq:lasota_0}), and using Assumption \ref{ass:pert} we have
\begin{equation}\label{eq:time-dep-lasota}
\|S_{X,Y, k}'\nu\|\leq 2\lambda^{-1}(1+\ve_1)^2\|\nu\|a_{z_k}a_{w_k}+BS_{X,Y, k}'\nu(1).
\end{equation}
Hence, for $\ve_1$ such that $2\lambda^{-1}(1+\ve_1)^2\leq \eta(1-\ve_1)^{-2}<1$ we can iterate the above inequality and obtain
\[
\begin{split}
&\| S_{X,Y, m-1}'\cdots  S_{X,Y, 0}'\nu\|
\leq \eta^m(1-\ve_1)^{2m}\|\nu\|\prod_{k=0}^{m-1}a_{z_k}a_{w_k}\\
&\quad+B\sum_{j=0}^{m-1}\eta^{j}(1-\ve_1)^{2j}\nu\left( S_{X,Y, 0}\cdots  S_{X,Y, m-1-j}1\right)\prod_{k=m-j}^{m-1}a_{z_k}a_{w_k}\\
&\leq\left[\eta^m\|\nu\|+(1-\eta)^{-1}B\right]\nu\left( S_{X,Y, 0}\cdots  S_{X,Y, m-1}1\right),
\end{split}
\]
which proves the Lemma with $\Cfive=1+(1-\eta)^{-1}B$.\footnote{Note that, for a probability measure, $1=\nu(1)=\nu(\partial_{\theta_i}\theta_i)\leq\|\nu\|$.}
\end{proof}

\section{Annealed Invariance Principle.}
\label{sec:invariance}

This section is devoted to proving an averaged invariance principle. 
This result is used in Section \ref{sec:clt-proof} to prove Lemma \ref{lm:hot}.

Consider the process
\begin{equation}\label{eq:path}
\hat X^N_t=\frac{1}{\sqrt{N}}\left\{\tilde X_{\lceil tN\rceil}+(tN-\lceil tN\rceil)\tilde\Delta_{\lceil tN\rceil}\right\}.
\end{equation}
Note that $\hat X^N_t\in\cC^0([0,1],\bR^d)$, by construction. In fact, Lemma \ref{lem:largedev} implies higher regularity.

\begin{lemma}\label{lem:holder}
The family of processes $\{\hat X^N\}\subset \cC^0([0,1],\bR^d)$ is tight.
\end{lemma}
\begin{proof}
Let $\varsigma\in (0,1/2)$,
\[
L_\varsigma(f):=\sup_{t,s\in[0,1]}\frac{\|f(t)-f(s)\|}{|t-s|^\varsigma},
\]
and $K_L^\varsigma:=\{f\in\cC^0([0,1],\bR^d)\;:\; f(0)=0,\; L_\varsigma(f)\leq L\}$.

By Lemma \ref{lem:largedev} it follows that, for each $N\in\bN$, $t\in [0,1]$ and $h\in [-t,1-t]$,
\begin{equation}\label{eq:holder}
\bP\left(\left\{\|\hat X^N_{t+h}-\hat X^N_t\|\geq L h^{\varsigma}\right\}\right)\leq e^{-CL^2h^{2\varsigma-1}}.
\end{equation}
In addition, if $\|\hat X^N_t\|+\|\hat X^N_{t+h}\|\leq L^{1-\varsigma}$, then the set in \eqref{eq:holder} is empty for all $h>L^{-1}$.
Now the result follows in complete analogy with the usual proof of the
H\"older continuity of the Brownian motion, based on applying the
above estimates to the dyadic rationals, yielding
\[
\bP\left(\left\{\hat X^N\not\in K_L^\varsigma\right\}\right)\leq e^{-CL}.
\]
Since $K_L^\varsigma$ are compact in $\cC^0([0,1],\bR^d)$ the tightness follows.
\end{proof}

Lemma \ref{lem:holder} also allows us to prove the invariance principle.
\begin{lemma}\label{lem:invariance}
For each probability measure $\nu\in\cB$ the process $\{\hat X^N_t\}$ converges in law to the Brownian motion with diffusion matrix $\Var$.
\end{lemma}
\begin{proof}
In view of Lemma \ref{lem:holder} we only need to check the convergence of finite dimensional distributions.
We consider two dimensional distributions, the general case being very similar. Accordingly, let
$t_1<t_2$ and fix $\xi_1, \xi_2.$ We have
$$ \bE_\nu\left(\exp(i\langle \xi_1, \hX_{t_1}^N\rangle+i\langle\xi_2 ,\hX_{t_2}^N\rangle)\right)=
\bE_\nu\left(\exp(i\langle[\xi_1+\xi_2],\hX_{t_1}^N\rangle)\exp(i\langle\xi_2 ,[\hX_{t_2}^N-\hX_{t_1}^N]\rangle)\right)=$$
$$ \bE_\nu\left(\bE\left(\exp(i\langle\xi_2, [\hX_{t_2}^N-\hX_{t_1}^N]\rangle)|\cF_{[t_1 N]}\right) 
\exp(i\langle[\xi_1+\xi_2], \hX_{t_1}^N\rangle)\right).$$
By Lemma \ref{lem:variation} and 
Proposition \ref{lem:clt-averaged} we 
have\footnote{In fact, Lemma \ref{lem:variation} considers two walks, 
yet the corresponding result for one walk can be obtained by integrating 
over the second walk. Moreover, for each smooth function $f:\bR^{2d}\to\bR$,
\[
\begin{split}
\bE_\nu(&f(\hat X^N_{t_2},\hat X^N_{t_1})\;|\;\cF_{[t_1N]})
=\bE_\nu\left(\bE\big(f(\hat X^N_{t_2},\hat X^N_{t_1})\;|\; \tX_{[t_1N]},\theta^{[t_1N]}\big)\;|\;\cF_{[t_1N]}\right)\\
&=\nu^X_{\tX_{[t_1N]},[t_1N]}\left(\bE(f(\hat X^N_{t_2},\hat X^N_{t_1}))\;|\; \tX_{[t_1N]},\theta^{[t_1N]})\right)\\
&=\bE_{X^N_{t_1},\nu^X_{\tX_{[t_1N]},[t_1N]}}\left(\bE(f(\hat X^N_{t_2-t_1},\hat X^N_{0})\;|\; \tX_{0},\theta^{0})\right)
=\bE_{X^N_{t_1},\nu^X_{\tX_{[t_1N]},[t_1N]}}\left(f(\hat X^N_{t_2-t_1},\hat X^N_{0})\right).
\end{split}
\]
Proposition \ref{lem:clt-averaged} can then be applied after 
translating $\nu^X_{\tX_{[t_1N]},[t_1N]}$ by $\tX_{[t_1N]}$.}
$$ \bE\left(\exp(i\langle\xi_2 ,[\hX_{t_2}^N-\hX_{t_1}^N]\rangle)|
\cF_{[t_1 N]}\right)=
\exp\left(-\frac{1}{2} \la  \xi_2,\Sigma^2 \xi_2\ra(t_2-t_1)(1+o(1))\right) $$
and so using Proposition \ref{lem:clt-averaged} again we obtain
$$ \bE_\nu\left(\exp(i\langle\xi_1, \hX_{t_1}^N\rangle+i\langle\xi_2, \hX_{t_2}^N\rangle)\right)\sim
e^{-\frac{1}{2} \left[\la  \xi_2,\Sigma^2 \xi_2\ra(t_2-t_1))+
\la  (\xi_1+\xi_2), \Sigma^2(\xi_1+\xi_2)\ra t_1\right]}$$
Thus, $(\hX_{t_1}^N, \hX_{t_2}^N)$ is asymptotically Gaussian with zero mean and
the variance predicted by the Brownian Motion. 
\end{proof}

\section{Proofs: Quenched CLT via the study of two random walks}
\label{sec:clt-proof}

The goal of this section is to establish Lemma \ref{lm:hot}.

Lemma \ref{lem:holder} shows that the distributions of the processes $(\hat X^N_t, \hat Y^N_t)$ are tight, 
hence they have accumulation points. Our next task is to characterize
such accumulation points. Let us consider any accumulation point
$(\hat X^\infty_{t},\hat Y^\infty_{t})$.  We will see that
$(\hat X_t^\infty,\hat Y^\infty_t)$ is a centered Gaussian random variables with variance
\begin{equation}\label{eq:vartwo}
\Vartwo:=t\begin{pmatrix}\Var&0\\0&\Var\end{pmatrix}.
\end{equation}
More precisely, if we define the second order differential operator 
$\Delta_{\Vartwo}:=\sum_{i,j}\Vartwo_{ij}\partial_i\partial_j$
we have the following.
\begin{prop}
\label{lem:heat}
For any $\psi\in\cC^{3}(\bR\times \bR^d\times \bR^d, \bR)$ we have
\[
\frac d{dt}\bE^2(\psi(t,\hat X^\infty_{t},\hat Y^\infty_{t}))=\bE^2(\partial_t\psi(t,\hat X^\infty_{t},\hat Y^\infty_{t})+
\frac 12\Delta_\Var\psi(t,\hat X^\infty_{t},\hat Y^\infty_{t}) ).
\]
More precisely, there exists $\beta\in(0,\frac 16)$ and $\vartheta\in (0,1-2\beta)$ such that, for all $N\in\bN$ and $t,h\in [0,1]$ such that $h>N^{\vartheta-1}$ we have
\[
\begin{split}
&\left|\bE^2\left(\psi(t+h,\hat X^N_{t+h},\hat Y^N_{t+h})-\psi(t,\hat X^N_{t},\hat Y^N_{t})
-h\left[\partial_t\psi(t,\hat X^N_{t},\hat Y^N_{t})+
\frac{1}{2} {\Delta_\Var\psi}(t,\hat X^N_{t},\hat Y^N_{t})\right] \right)\right|\\
&\quad \leq C\|\psi\|_{\cC^3}(N^{-\beta}h+h^{\frac 32}+N^{-\frac 12}).
\end{split}
\]
\end{prop}
Thanks to the above Proposition for each $\phi\in\cC_0^3(\bR^d\times \bR^d)$ we can define $\psi$ by
\begin{equation}\label{eq:backheat}
\begin{split}
&\partial_t\psi+\frac 12\Delta_\Var\psi=0\\
&\psi(1,x,y)=\phi(x,y)
\end{split}
\end{equation}
and, by applying Proposition \ref{lem:heat} with the choice $h=\lceil N^{2\beta}\rceil^{-1}$, obtain the wanted result:
\[
\begin{split}
\bE^2&(\phi(\hat X^N_{1},\hat Y^N_{1}))=\bE^2(\psi(1,\hat X^N_{1},\hat Y^N_{1}))\\
&=\bE^2(\psi(0,\hat X^N_{0},\hat Y^N_{0}))+\sum_{i=0}^{h^{-1}-1}\bE^2(\psi((i+1)h,\,\hat X^N_{(i+1)h},\hat Y^N_{(i+1)h})-\psi(ih,\,\hat X^N_{ih},\hat Y^N_{ih}))\\
&=\bE^2(\psi(0,0,0))+\cO(\|\psi\|_{\cC^3}N^{-\beta})=\bE_{\cN(0,\Vartwo)}(\phi)+\cO(\|\phi\|_{\cC^3}(N^{-\beta}+N^{-\frac 12+2\beta})),
\end{split}
\]
where we have used the explicit solution of \eqref{eq:backheat}.\footnote{Indeed, equation \eqref{eq:backheat}  is just the backward heat equation, thus for $t\in (0,1)$
\[
\psi(t,x,y)=\frac1{(4\pi)^{d}\det(\Sigma_2^{-1}) (1-t)^d}\int_{\bR^{2d}} e^{-\frac{\langle(x-z,y-w),\Sigma_2^{-2}(x-z,y-w)\rangle}{2(1-t)}}\phi(z,w) \;dz\,dw.
\]
}
Remembering the form of $\Vartwo$ (see \eqref{eq:vartwo}), 
Lemma \ref{lm:hot}, 
and hence Theorem \ref{thm:main}, follow.

\begin{proof}[{\bfseries Proof of Proposition \ref{lem:heat}}]
We start by the following Taylor expansion
\begin{equation}\label{eq:hoihoi}
\begin{split}
\bE^2&(\psi(t+h,\hat X^N_{t+h},\hat Y^N_{t+h}))-\bE^2(\psi(t,\hat X^N_{t},\hat Y^N_{t}))=\bE^2(\partial_t\psi(t,\hat X^N_{t},\hat Y^N_{t}))h\\
&+\bE^2(\partial_x\psi(t,\hat X^N_{t},\hat Y^N_{t})\cdot(\hat X^N_{t+h}-\hat X^N_t)
+\partial_y\psi(t,\hat X^N_{t},\hat Y^N_{t})\cdot(\hat Y^N_{t+h}-\hat Y^N_t))\\
&+\frac 12\bE^2((\hat X^N_{t+h}-\hat X^N_t)\cdot\partial_x^2\psi(t,\hat X^N_{t},\hat Y^N_{t})\cdot(\hat X^N_{t+h}-\hat X^N_t))\\
&+\bE^2((\hat Y^N_{t+h}-\hat Y^N_t)\cdot\partial_{xy}\psi(t,\hat X^N_{t},\hat Y^N_{t})\cdot(\hat X^N_{t+h}-\hat X^N_t))\\
&+\frac 12\bE^2((\hat Y^N_{t+h}-\hat Y^N_t)\cdot\partial_y^2\psi(t,\hat X^N_{t},\hat Y^N_{t})\cdot(\hat Y^N_{t+h}-\hat Y^N_t))\\
&+\frac 12\bE^2(\partial_{tx}\psi(t,\hat X^N_{t},\hat Y^N_{t})\cdot(\hat X^N_{t+h}-\hat X^N_t)) h\\
&+\frac 12\bE^2(\partial_{ty}\psi(t,\hat X^N_{t},\hat Y^N_{t})\cdot(\hat Y^N_{t+h}-\hat Y^N_t)) h\\
&+\cO\left(|\partial_t^2\psi|_\infty h^2+ \|\psi\|_{\cC^3}\left[\bE(\|\hat X^N_{t+h}-\hat X^N_t\|^3)+\bE(\|\hat Y^N_{t+h}-\hat Y^N_t\|^3)\right]\right).
\end{split}
\end{equation}
Next, we will analyze the terms in equation \eqref{eq:hoihoi} one by one. 

First of all note that the remainders are of order $h^{\frac 32}$.\footnote{\label{foot:momenta} In fact Lemma \ref{lem:largedev} implies, for each $p\in\bN$, that 
\[
\begin{split}
\bE(\|\hat X^N_{t+h}-\hat X^N_t\|^p)&\leq N^{\frac p2}h^p\bE\left(\left\|\frac{1}{\lceil hN\rceil}(\tilde X_{\lceil (t+h)N\rceil}-\tilde X_{\lceil tN\rceil})\right\|^p\right)\\
&\leq Ch^{\frac p2}+C_pN^{\frac p2}h^p\int_{(Nh)^{-\frac 12}}^\infty x^{p-1}\bP\left(\left\{ 
\left\|\frac{1}{\lceil hN\rceil}(\tilde X_{\lceil (t+h)N\rceil}-\tilde X_{\lceil tN\rceil})\right\|\geq x\right\}\right)dx\\
&\leq Ch^{\frac p2}+C_pN^{\frac p2}h^p\int_{(Nh)^{-\frac 12}}^\infty x^{p-1} e^{-Cx^2hN} dx\leq C_p h^{\frac p2}.
\end{split}
\]
}
Let us start the estimates with the term $\bE(\partial_x\psi(t,\hat X^N_{t},\hat Y^N_{t})\cdot(\hat X^N_{t+h}-\hat X^N_t))$. 
To this end we consider the $\sigma$-algebra $\cF^{XY}_m$ generated by $\{X_1,\dots, X_m,Y_1,\dots,Y_m\}$.
Setting $\tg=g-v$, $\ell_t:=\lceil tN\rceil+1$ we can write
\[
\begin{split}
\bE((\hat X^N_{t+h}-\hat X^N_t)\;|\; \cF^{XY}_{\ell_{t}})&=N^{-\frac 12}\sum_{k=\ell_t}^{\ell_{t+h}}\bE(\tilde \Delta_k^X\;|\;\cF^{XY}_{\ell_{t}})+\cO(N^{-\frac 12})\\
&=N^{-\frac 12}\sum_{k=\ell_t}^{\ell_{t+h}}
\bE( (S^{k-\ell_{t}}\tg)\circ \tau^{X_{\ell_{t}}}\;|\;\cF^{XY}_{\ell_{t}})+\cO(N^{-\frac 12})\\
&=N^{-\frac 12}\sum_{k=\ell_t}^{\ell_{t+h}}
[(S^{k-\ell_{t}}\tau^{X_{\ell_{t}}})'(\mue)^{XY}_{\ell_t}](\tg)+\cO(N^{-\frac 12}).
\end{split}
\]
where, by Lemma \ref{lem:variation},
$\|(\mue)^{XY}_{\ell_t} \|\leq C$ and hence $\|(\tau^{X_{\ell_{t}}})'(\mue)^{XY}_{\ell_t} \|\leq C.$
From Theorem \ref{thm:mixing} and the fact that $\muinv(\tg)=0$ it follows that 
\begin{equation}\label{eq:first_term}
|\bE^2(\partial_x\psi(t,\hat X^N_{t},\hat Y^N_{t})\cdot(\hat X^N_{t+h}-\hat X^N_t))|\leq CN^{-\frac 12}\|\psi\|_{\cC^1}.
\end{equation}
In complete analogy we have
\begin{equation}\label{eq:second_term}
|\bE^2(\partial_y\psi(t,\hat X^N_{t},\hat Y^N_{t})\cdot(\hat Y^N_{t+h}-\hat Y^N_t))|\leq CN^{-\frac 12}\|\psi\|_{\cC^1}.
\end{equation}
The quadratic terms involving 
$\partial_{tx}$ and $\partial_{ty}$ 
are estimated in the same manner yielding terms of order
$N^{-\frac 12} h\|\psi\|_{\cC^2}$. 
The quadratic terms involving only $X$ or only $Y$ 
yield the following 
\begin{equation}\label{eq:third_term}
\begin{split}
&\bE^2((\hat X^N_{t+h}-\hat X^N_t)_i
(\hat X^N_{t+h}-\hat X^N_{t})_j\partial_{x_ix_j}\psi)
=(\Var)_{ij}\bE^2(\partial_{x_ix_j}\psi)h+\cO\left(\frac{\|\psi\|_{\cC^2}}{\sqrt
N}\right)\\
&\bE^2((\hat Y^N_{t+h}-\hat Y^N_t)_i (\hat Y^N_{t+h}-\hat
Y^N_{t})_j\partial_{y_iy_j}\psi)
=(\Var)_{ij}\bE^2(\partial_{y_iy_j}\psi)h+\cO\left(\frac{\|\psi\|_{\cC^2}}{\sqrt
N}\right).
\end{split}
\end{equation}
Indeed, in analogy with what we have done before, we can condition with respect to $\cF^{XY}_{\ell_t}$ and
\begin{equation}\label{eq:details2}
\begin{split}
\bE^2&((\hat X^N_{t+h}-\hat X^N_t)_i
(\hat X^N_{t+h}-\hat X^N_{t})_j\;|\;\cF^{XY}_{\ell_t})=\frac 1N\left\{ \sum_{k=\ell_t}^{\ell_{t+h}-1}\bE^2((\tDelta^X_k)_i(\tDelta^X_k)_j\;|\;\cF^{XY}_{\ell_t})\right.\\
&\left. +\sum_{k=\ell_t}^{\ell_{t+h}-1}\sum_{l=\ell_t}^{k-1}\bE^2((\tDelta^X_k)_i(\tDelta^X_l)_j+(\tDelta^X_l)_i(\tDelta^X_k)_j\;|\;\cF^{XY}_{\ell_t})\right\}\\
&+N^{-1}\cO\left(1+\sum_{l=\ell_t+1}^{\ell_{t+h}-1}\bE^2((\tDelta^X_{\ell_{t+h}-1})_i(\tDelta^X_l)_j+(\tDelta^X_{\ell_{t+h}-1})_j(\tDelta^X_l)_i\;|\;\cF^{XY}_{\ell_t})\right),
\end{split}
\end{equation}
where the boundary terms of the type $\bE^2((\tDelta^X_{\ell_{t}-1})_i(\tDelta^X_l)_j\;|\;\cF^{XY}_{\ell_t})$ have been estimated as in \eqref{eq:first_term}.
To estimate such an expression note that
\[
\bE^2((\tDelta^X_k)_i(\tDelta^X_l)_j\;|\;\cF^{XY}_{\ell_t})=
(S^{j-\ell_t}\tau^{X_{\ell_t}})'\nu^{XY}_{\ell_t}\left(\tg_jS^{k-j}\tg_k\right)
\]
Since Lemma \ref{lem:variation} and Theorem \ref{thm:mixing} imply $\|(S^{j-\ell_t}\tau^{X_{\ell_t}})'\nu^{XY}_{\ell_t}\|\leq C\|\mue\|$ and also that $\|\tg_j(S^{j-\ell_t}\tau^{X_{\ell_t}})'\nu^{XY}_{\ell_t}\|\leq C\|\mue\|$ we can apply Theorem \ref{thm:mixing} 
\[
\begin{split}
\left|\bE^2((\tDelta^X_k)_i(\tDelta^X_l)_j\;|\;\cF^{XY}_{\ell_t})\right|= &\nu^{XY}_{\ell_t}(S^{j-\ell_t}\tau^{X_{\ell_t}}\tg_j)\muinv(\tg_k)+\cO(\eta^{k-j})=\cO(\eta^{k-j})\\
\left|\bE^2((\tDelta^X_k)_i(\tDelta^X_l)_j\;|\;\cF^{XY}_{\ell_t})\right|=&\nu^{XY}_{\ell_t}(1)\muinv(\tg_jS^{k-j}\tg_k)+\cO(\eta^{j-\ell_t})\\
=&\muinv(\tg_jS^{k-j}\tg_k)+\cO(\eta^{j-\ell_t}).
\end{split}
\]
Using such estimates in \eqref{eq:details2} and remembering formula \eqref{eq:sec-der}  (for $\ddot\alpha_0=\Var$) the first equation of 
\eqref{eq:third_term} follows. The second equation
of \eqref{eq:third_term} is proven in complete analogy.

Finally, we must deal with the mixed quadratic 
term. 
\[
\bE^2((\hat X^N_{t+h}-\hat X^N_t)_i (\hat Y^N_{t+h}-\hat Y^N_{t})_j)\;|\;\cF_{\ell_t}^{XY})
=\!\!\!\!\sum_{k,m=\ell_t-1}^{\ell_{t+h}}\!\!\!\!\frac{\bE^2((\tilde \Delta_k^X)_i
(\tilde\Delta_m^Y)_j \;|\;\cF^{XY}_{\ell_{t}})}N+\cO( \frac 1{\sqrt N}).
\]
If $|k-m|\geq A\ln N$, then, for $A>\ln \eta^{-1}$, by Lemma \ref{lem:variation} and Theorem \ref{thm:mixing} it follows that 
\[
\begin{split}
\bE^2((\hat X^N_{t+h}-\hat X^N_t)_i (\hat Y^N_{t+h}-\hat Y^N_{t})_j)\;|\;\cF_{\ell_t}^{XY})
=&\!\!\!\!\sum_{\substack{\ell_t+2A\ln N\leq k\leq\ell_{t+h}\\|m-k|\leq
A\ln N}}\!\!\!\!\!\!\!\!\!\! \frac{\bE^2((\tilde \Delta_k^X)_i
(\tilde\Delta_m^Y)_j \;|\;\cF^{XY}_{\ell_{t}})}{N}\\
&\quad+\cO( N^{-\frac 12}). 
\end{split}
\]
Next, suppose that $|m-k|\leq A\ln N$ and $|X_{k-A\ln N}- Y_{k-A\ln N}|> 4\Czero A\ln N$.

Assume, to fix our ideas, that $k\leq m$. Then, for all times 
$l$ such that $|l-k|\leq A\ln N$, the two walks explore disjoint parts of the environment. Thus, we can consider the process started at time $k-A\ln N$ with the conditional measure $(\mue)^{XY}_{k-A\ln N}$ and with the walks starting from $a=X_{k-A\ln N}, b=Y_{k-A\ln N}$, $\|a-b\|>4\Czero A\ln N$.
If we set ${\bf f}=S^{A\ln N}\tau^{X_k}\tg_i$ and ${\bf h}=S^{m-k+A\ln N}\tau^{Y_k}\tg_j$ we have that the two functions depend on different sets of variables (let $B\subset \bZ^d$ be the set of variables on which $\bf f$ depends and $B'$ the ones relative to $\bf h$) and
\[
\bE^2((\tilde \Delta_k^X)_i (\tilde\Delta_m^Y)_j\;|\;\cF^{XY}_{k-A\ln N})=(\mue)^{XY}_{k-A\ln N}({\bf f}\,{\bf h}).
\]
We can then define its Newtonian potential of $\bf f$,
\begin{equation}\label{eq:newton}
\begin{split}
\Psi(\theta)&=\frac{1}{|B|(|B|-2)\alpha_{|B|}}\int_{I^{\bZ^d}}\|\theta^B-\vartheta^B\|^{-|B|+2}{\bf f}(\vartheta)d\vartheta^B \otimes_{j\not\in B}\mu_0(d\vartheta_j)\\
&=\frac{1}{|B|(|B|-2)\alpha_{|B|}}\int_{I^{B}}\|\theta^B-\vartheta^B\|^{-|B|+2}{\bf f}(\vartheta^B)d\vartheta^B.
\end{split}
\end{equation}
where $\theta^B=(\theta_l)_{l\in B}$ and $\alpha_{l}$ is the volume of the unit ball in $\bR^l$. It is well known that, for $\theta^B$ in the interior of $I^B$
\[
\sum_{l\in B}\partial_{\theta_l\theta_l}\Psi=\bf f.
\]
Thus, remembering Lemma \ref{lem:variation}, we can write\footnote{Remember that the marginal of $(\mue)^{XY}_\ell$ on $B\cup B'$ is absolutely continuous with respect to Lebesgue, hence the boundary of $I^B$ has zero measure, moreover $\partial_{\theta_l}\Psi$ is a continuous function on $I^{\bZ^d}$.}
\[
\begin{split}
&\left|\bE^2((\tilde \Delta_k^X)_i (\tilde\Delta_m^Y)_j\;|\;\cF^{XY}_{k-A\ln N})\right|\leq\sum_{l\in B}\left|(\mue)^{XY}_{k-A\ln N}(\partial_{\theta_l}(\partial_{\theta_l}\Psi\cdot {\bf h}))\right|\\
&\quad\leq |B|\cdot \|(\mue)^{XY}_{k-A\ln N}\|\cdot \sup_{l\in B} |\partial_{\theta_l}\Psi\cdot {\bf h}|_\infty
\leq CA^d\Czero^d|\ln N|^d\cdot |g|_\infty\cdot\sup_{l\in B} |\partial_{\theta_l}\Psi|_\infty.
\end{split}
\]
By \eqref{eq:newton} we have, for $l\in B$,
\[
\partial_{\theta_l}\Psi(\theta)=\int_{I^{\bZ^d}}\frac{\theta_l-\vartheta_l}{|B|(|B|-2)\alpha_{|B|}\cdot\|\theta^{B}-\vartheta^B\|^{|B|}} \;{\bf f}(\vartheta)d\vartheta^B \otimes_{j\not\in B}\mu_0(d\vartheta_j)=:\nu^\theta_l({\bf f}).
\]
Unfortunately, $\nu^\theta_l\not \in\cB$ due to the singularity of the kernel. To take care of this problem we need to isolate the singularity. For each $r>0$ let $\chi_r\in\cC^\infty(\bR^B,[0,1])$ such that $\chi_r(\theta^B)=0$ for all $\|\theta^B\|\leq r$ and $\chi_r(\theta^B)=1$ for all $\|\theta^B\|\geq2r$. Clearly $\chi_r$ can be chosen radial and so that $\sup_l|\partial_{\theta_l}\chi_r|_\infty\leq Cr^{-1}$.
We then define 
\[
\nu^\theta_{l,r}(\phi):=\int_{I^{\bZ^d}}\frac{(\theta_l-\vartheta_l)\cdot\chi_r(\theta-\vartheta)}{|B|(|B|-2)\alpha_{|B|}\cdot\|\theta^{B}-\vartheta^B\|^{|B|}} \;\phi(\vartheta)\;\;d\vartheta^B \otimes_{j\not\in B}\mu_0(d\vartheta_j)
\]
 and $\mu^\theta_{{l,r}}(f):=\nu^\theta_l(\phi)-\nu^\theta_{l,r}(\phi)$. A direct computation shows that $|\mu^\theta_{l,r}|\leq C r(A\Czero\ln N)^{-d}$ and $\|\nu^\theta_{l,r}\|\leq C(A\Czero\ln N)^{-d}\ln r^{-1}$. Since $|{\bf f}|\leq 2|g|_\infty$ we can finally use Theorem \ref{thm:mixing} to estimate
\[
\begin{split}
&\left|\bE^2((\tilde \Delta_k^X)_i (\tilde\Delta_m^Y)_j\;|\;\cF^{XY}_{k-A\ln N})\right|
\leq C\left\{ r+\sup_{\substack{l\in B\\\theta\in I^B}}|(\tau^{X_k})'\nu^\theta_{l,r}(S^{A\ln N}\tg_i)|\right\}\\
&\quad \leq C\left\{ r+\muinv(\tg_i)+\sup_{\substack{l\in B\\\theta\in I^B}}\|
(\tau^{X_k})'\nu^\theta_{l,r}\|\eta^{A\ln N}\right\}\\
&\quad\leq C\left\{ r+N^{-1}\ln r^{-1}\right\}\leq CN^{-1}\ln N \,,
\end{split}
\]
where we have chosen $r=N^{-1}$.

In conclusion,
\begin{equation}\label{eq:fourth_term}
\begin{split}
\big|\bE^2(&(\hat X^N_{t+h}-\hat X^N_t)_i \partial_{x_iy_j}\psi(t,\hat X^N_t,\hat Y^N_t)(\hat Y^N_{t+h}-\hat Y^N_{t})_j)\big|\leq CN^{-\frac 12}\|\psi\|_{\cC^2}\\
&+ A\ln N\;N^{-1}\bE^2(\Card\{t\leq N: \|X_t-Y_t\|\leq 4\Czero A\ln N\})\|\psi\|_{\cC^2}.
\end{split}
\end{equation}
In Section \ref{sec:twowalks} we prove the following bound.
\begin{lemma}
\label{lem:NL} 
Let $A$ be a large constant and set $L_N:=A\ln N$. 
There exists  $\delta_0\in (0,1)$ such that
\begin{equation}
  \label{NL}
  \bE^2(\Card\{t\leq N: \|X_t-Y_t\|\leq L_N\})\leq C N^{\delta_0} \quad (N\in\bN).
\end{equation}
\end{lemma}
Lemma \ref{lem:NL} allows to estimate the last term in the right hand side of \eqref{eq:fourth_term} by
$$ C N^{\delta_0-1}\ln N \|\psi\|_{\cC^2}=C h (h^{-1} N^{\delta_0-1}) \ln N \|\psi\|_{\cC^2}, $$
proving the proposition by choosing $\beta=\frac{1-\delta_0}6$ and $\vartheta=1-3\beta$.
\end{proof}

\section{Two walks estimates}
\label{sec:twowalks}
In Section \ref{sec:clt-proof} we proved that Lemma \ref{lm:hot} 
(and hence Theorem \ref{thm:main}) holds provided the average number of times two walks come closer than $A\ln N$ in time $N$ is smaller than $N^{\delta_0}$ for some $\delta_0\in(0,1)$. The purpose of this section is to prove such an estimate and therefore conclude the argument.

\subsection{On the number of close encounters}\label{sec:two-crossing}
The proof of inequality \eqref{NL} can be reduced to the following simpler inequality.
\begin{lemma}
  \label{lem:Ex}
  There exist $\rho\in(0,1), \Cfour>0$ such that 
  for any $m\in\bN$ and for any $a, b$ such that
  $\|a-b\|> L_N$, we have 
  \begin{equation}
    \label{Ex}
    \bP^2\left(\left\{\|X_j-Y_j\|>L_N \;\;\;\forall\, j\in\{m,\dots,m+N\}\right\}\;|\; 
      X_m=a, Y_m=b\right)\geq \frac{\Cfour}{N^{\rho}} ,
  \end{equation}
  (Here $\bP^2$ is the underlying probability
  for the process $(\theta_t,X_t,Y_t)$ started with $\theta_0$ distributed according to $\mue$).
\end{lemma}
 We postpone the proof of the above Lemma until finishing the proof of \eqref{NL}.
 
 \begin{proof}[\bfseries Proof of Lemma \ref{lem:NL}]
 Notice that Assumption \ref{ass:ellipticity} implies that the walks can move in different directions with positive probability. In particular, there exists $\gamma>0$ such that for each
$a,b\in\bZ^d$, $m\in\bN$ and $\delta>0$,
  \begin{equation}
    \label{eq:fast-escape}
    \bP^2\left(\left\{\|X_{m+L_N^2}-Y_{m+L_N^2}\|\geq
        L_N\right\}\;\bigg|\; X_m=a,Y_m=b\right)\geq c(\delta)\gamma^{\delta L_N},
  \end{equation}
  Indeed 
$$   
\bP^2\left(\left\{\|X_{m+\delta L_N}-Y_{m+\delta L_N}\|\geq
        \delta L_N\right\}\;\bigg|\; X_m=a,Y_m=b\right)\geq
\gamma^{\delta L_N} 
$$
  the latter being the probability of one fixed path in which $X_i,Y_i$ get
  further and further apart at each step.
On the other hand 
\begin{equation}
\label{LNCLT}
     \bP^2\left(\left\{\|X_{m+L_N^2-\delta L_N}-Y_{m+L_N^2-\delta L_N}\|\geq
        L_N\right\}\;\bigg|\; \|X_m-Y_m\|\geq \delta L_N \right)\geq c(\delta) , 
 \end{equation}  
 To verify \eqref{LNCLT}  let $W^{(1)}(t)$ and $W^{(2)}(t)$ be independent Brownian Motions
 such that $W^{(2)}(0)-W^{(1)}(0)=\bv,$ where $\|\bv\|=\delta.$ Let $\hat\bv:=\bv\,\|\bv\|^{-1}$, and
 \[
 c(\delta):=\frac{1}{2}P\left(\|W^{(2)}(1)-W^{(1)}(1)\|>1, \quad \text{and for all } t\in [0,1] \right.
 \]
 \[
 \la W^{(2)}(t), \hat\bv \ra>     \la   W^{(2)}(0), \hat\bv \ra-\frac{\delta}{3}     
  \text{ and } \left. \la W^{(1)}(t), \hat\bv \ra<      \la W^{(1)}(0), \hat\bv \ra-\frac{\delta}{3}\right)
  \]
Observe that the invariance principle established in Lemma
\ref{lem:invariance}, Lemma \ref{lem:variation} and the fact that
local  dynamics are independent implies a two particle invariance
principle as long as the walkers explore disjoint regions in the phase space.
Therefore  the probability that two walkers grow $L_N$ apart exploring disjoint regions of the phase space is at least $c(\delta)$ for large $N$ proving \eqref{LNCLT}.  
Choose $\varrho<1-\rho$. Then,
  \begin{equation}
    \label{eq:long-escape} 
    \begin{split}
      \bP^2&\left(\left\{\sup_{m\leq i\leq m+N^\varrho}\|X_i-Y_i\|\leq
          L_N\right\}\;\bigg|\; X_m=a,Y_m=b\right)\\
          &\leq \prod_{j=1}^{N^\varrho L_N^{-2}}(1-c(\delta)\gamma^{\delta L_N})
      \leq e^{-c(\delta) \gamma^{\delta L_N}L_N^{-2}N^\varrho}\leq
            e^{-CN^{\varrho/2}},
    \end{split}
  \end{equation}
provided that $\delta$ is sufficiently small.

  Next, consider the sets $B_R^-:=\{(x,y): \|x-y\|\leq R\}$,
  $B_R^+:=\{\|x-y\|> R\}$ and the stopping times, for $k>0$,
  \[
  \begin{split}
    s_0&:=\inf\left\{j\in \bN\;:j>0,\;\;(X_{j-1},Y_{j-1})\in B^-_{L_N},
      (X_{j},Y_{j})\in B^+_{L_N}\right\},\\
    s_{2k}&:=\inf\left\{j\in \bN\;:\; j> s_{2k-2},\; (X_{j-1},Y_{j-1})\in B^-_{L_N},
      (X_{j},Y_{j})\in B^+_{L_N}\right\},\\
    s_1&:=\inf\left\{j\in \bN\;:\;j>s_0,\;(X_{j-1},Y_{j-1})\in B^+_{L_N},\; 
      (X_{j},Y_{j})\in B^-_{L_N}\right\},\\ 
    s_{2k+1}&:=\inf \left\{j\in\bN\;:\; j> s_{2k-1},\, (X_{j-1},Y_{j-1})\in
      B^+_{L_N},\; (X_{j},Y_{j})\in B^-_{L_N}\right\}.
  \end{split}
  \]
  Clearly, $s_{2k}<s_{2k+1}<s_{2k+2}$ and $s_k>k$. As $X_0=Y_0$, these
  stopping times are adapted to the filtration $\cF^{XY}_t$. Note that the $s_{2k}$ are upcrossing times hence $\|X_t-Y_t\|\leq L_N$ for all the $t\in\{s_{2k-1}, \dots, s_{2k}-1\}$.
  With this notation, \eqref{eq:long-escape} implies
  \begin{displaymath}
    \bP^2\left(\left\{\sup_{i\leq N}
        (s_{2i}-s_{2i-1})>N^\varrho\right\}\right)
    \leq
    N\sup_{i\leq N}\bP^2 \left(\left\{s_{2i}-s_{2i-1}>N^\varrho\right\}\right)
    \leq
    Ne^{-N^{\varrho/2}}.
  \end{displaymath}
  Let us set $J:=\inf \{k\in\bN\;:\;s_k\geq N\}$, clearly $J\leq N$,
  \begin{equation}
    \label{eq:NL0}\bE^2(\Card\{n<N\;:\;\|X_n-Y_n\|\leq L_N\})
    \leq 
    N^2e^{-N^{\varrho/2}}+N^\varrho\bE^2(J).
  \end{equation}
  
 It remains to investigate the length of the intervals of time in which the
  two walks are further apart than $L_N$.  Let $S_n:=\{\sup_{k\leq
    n}(s_{2k+1}-s_{2k})<N\}$, and denote by $\cF_{s_{2k}}^{XY}$ the
  $\sigma$-algebra associated to the filtration $\cF_t^{XY}$ and the stopping time
  $s_{2k}$. Then, by \eqref{Ex},
  \begin{displaymath}
    \begin{split}
      \bP^2(\{J>n\})
      &\leq
      \bE^2\left(\Id_{S_{n}}\right)
      =
      \bE^2(\Id_{S_{n-1}}\bP^2
      (\{s_{2n+1}-s_{2n}<N\}\;|\;S_{n-1}))\\
      &=
      \bE^2\left(\Id_{S_{n-1}}\bP^2
        \left(\{s_{2n+1}-s_{2n}<N\}\;|\;
          \cF^{XY}_{s_{2n}}\right)\right)\\
      &\leq 
      \left(1-\frac {C_0}{N^\rho}\right)\bE^2(\Id_{S_{n-1}})
      \leq \dots \leq
      \left(1-\frac {C_0}{N^\rho}\right)^{n}.
    \end{split}
  \end{displaymath}
  Thus, letting $1-\varrho>\alpha>\rho$, it follows that
  \[
  \bP^2(\{J>N^\alpha\})\leq Ce^{-C_0N^{\alpha-\rho}}.
  \]
  which means that $\bE^2(J)\leq
  N^\alpha+N\,\bP^2(\{J>N^\alpha\})\leq CN^\alpha$. In view of
  \eqref{eq:NL0} this proves \eqref{NL} provided we have chosen $\delta_0$ so
  that $\varrho+\alpha<\delta_0$.
\end{proof}

Our program is thus completed once we prove \eqref{Ex}. To this end an
intermediate result is needed.
\begin{lemma}
  \label{LmFarInd}
Given $R>0$, take two points $a_R$ and $b_R$ such that $\|a_R-b_R\|=R.$ 
Consider two walks starting at $a_R$ and $b_R$ respectively with the environment given by the probability distribution $\nu\in\cB$ and define the stopping time $\tau_{\delta,R}$ as the first time $n>0$ such that
  \begin{displaymath}
    \|X_n-Y_n\|\leq\frac{R}{1+\delta} \text{ or }
    \|X_n-Y_n\|\geq(1+\delta)R.     
  \end{displaymath}
For each $\Cnine>0$, if $\|\nu\|\leq \Cnine$,  then there exist $R_\delta\in\bR_+,$ $c_1, c_2>0$ such that for each $R\geq R_\delta$
  \begin{displaymath}
        \bP_\nu^2\left(\{\|X_{\tau_{\delta,R}}-Y_{\tau_{\delta,R}}\|\geq (1+\delta)R\}\right)\geq
      \frac{1}{2+\delta}-c_1 e^{-c_2/\delta}  \,.
    \end{displaymath}
\end{lemma}
\begin{proof}
Let $T_R$ be the first time $n>0$ such that 
$$ \max(\|X_n-X_m\|, \|Y_n-Y_m\|)\geq \frac{R}{2}. $$
Then, by Section \ref{sec:invariance} the pair
$$ \left(\frac{X_{\min(T_R, t R^2)}-X_m}{R}, \frac{Y_{\min(T_R, t R^2)}-X_m}{R} \right) $$
is asymptotic, when $R\to\infty$, to a pair of independent Brownian Motions 
$$(W^{(1)}(t), W^{(2)}(t))\in\bR^{d}\times \bR^d\quad W^{(1)}(0)=0, \quad 
\|W^{(2)}(0)\|=1 $$ 
stopped at time $T$ when 
one of them wanders more than $1/2$ from its starting position.\footnote{More precisely, the fact that each component is approximately Brownian comes from
Section \ref{sec:invariance} while independence is due to the fact that the walkers 
explore non-intersecting regions in the phase space and can be proven by the same arguments use to estimate \eqref{eq:fourth_term}.}
Let $\btau$ be the first time $\|W^{(1)}(t)-W^{(2)}(t)\|=(1+\delta)$ or
$\|W^{(1)}(t)-W^{(2)}(t)\|=(1+\delta)^{-1}.$ 
Recall that
$$ 
P(\{\|W^{(1)}(\btau)-W^{(2)}(\btau)\|=(1+\delta)\})\geq\frac{1}{2+\delta} 
$$
(the worst case is then $d=1$, see e.g. \cite{RY}, Section XI.1). On the other hand
\begin{equation}
\label{NearFar} 
 P(T<\btau)\leq P(\btau>\delta)+P(T<\delta) 
\end{equation} 
and both terms are $O(e^{-c/\delta}),$ the first one because
$$P(\btau>(k+1)\delta^2|\btau>k\delta^2)\leq\gamma<1$$ 
and the second one by Hoeffding's inequality
(see e.g. \cite{GS}). 
Now\footnote{Here $c_1$ should be taken a little bit larger than the implied constants
in \eqref{NearFar} to take into account that our process is only approximated by
Brownian Motion.}
\[
\begin{split}
 &\bP(\{\|X_{\tau_{\delta, R}}-Y_{\tau_{\delta, R}}=R(1+\delta)\})\\
 &\quad\geq \bP(\{\|X_{\tau_{\delta, R}}-Y_{\tau_{\delta, R}}=R(1+\delta)\text{ and }
\tau_{\delta, R}<T_R\})\geq
\frac{1}{2+\delta}-c_1 e^{-c_2/\delta} . 
\end{split}
\] 
\end{proof}

We now use the following comparison criterion (proved in section \ref{sec:comparison}).

\begin{lemma}
\label{PrComp}
Suppose $\xi_1, \xi_2\dots \xi_n\dots$ is a random process such that $\xi_n=\pm 1 $
and for all $n$ 
$$ P(\xi_n=1 | \xi_1\dots \xi_{n-1})\geq p . $$
Let $\txi_1, \txi_2 \dots \txi_n\dots $ be iid random variables such that $\txi_n=\pm 1,$
and $P(\txi_n=1)=p.$ Let
$$ \cX_n=\sum_{j=1}^n \xi_n+  \cX_0\quad  \tcX_n=\sum_{j=1}^n \txi_n+ \tcX_0. $$  
Then for any $ \alpha_1<\alpha< \alpha_2$ 
$$ P(\cX_k\text{ reaches } \alpha_2 \text{ before } \alpha_1|\cX_0=\alpha)\geq
P(\tcX_k\text{ reaches } \alpha_2 \text{ before } \alpha_1|\tcX_0=\alpha) $$
\end{lemma}
Recall that by Gambler's Ruin Formula for $p\neq 1/2$ 
\begin{equation}
\label{GR}
P(\tcX_k\text{ reaches } \alpha_2 \text{ before } \alpha_1|\tcX_0=\alpha)=
\frac{\left(\frac{p}{1-p}\right)^{\alpha_2-\alpha}-\left(\frac{p}{1-p}\right)^{\alpha_2-\alpha_1}}
{1-\left(\frac{p}{1-p}\right)^{\alpha_2-\alpha_1}} 
\end{equation}

\begin{proof}[\bfseries Proof of Lemma \ref{lem:Ex}]
Let $X_m=a$ and $Y_m=b$ with
$\|a-b\|\geq L_N$ and $\kappa\in(\frac 12,1)$.

Using ellipticity of Assumption \ref{ass:ellipticity} for the first $\delta L_N$ steps we see that
with probability greater than $N^{-c\delta}$ our walkers move distance $(1+\delta) L_N$ apart 
without getting within distance $L_N$ from each other. 
Let $\tau_1$ be the first time after $m$ when our walkers move distance $(1+\delta) L_N$ apart
and let $\tau_{n+1}$ be the first time after $\tau_n$ when 
$$ \|X_j-Y_j\|\geq (1+\delta) \|X_{\tau_n}-Y_{\tau_n}\| \text{ or } 
\|X_j-Y_j\|\leq (1+\delta)^{-1} \|X_{\tau_n}-Y_{\tau_n}\|.$$ 
Applying Lemma \ref{PrComp} to $\cX_n=\frac{\ln\|X_{\tau_n}-Y_{\tau_n}\|}{\ln(1+\delta)}$ with $\alpha_1=\frac{\ln L_N}{\ln(1+\delta)}$, $\alpha=\alpha_1+1$, $\alpha_2=\frac{\kappa\ln N}{\ln(1+\delta)}$ and
using Lemma \ref{LmFarInd}, taking into account Lemma \ref{lem:variation}, to estimate the probability of moving apart we conclude from \eqref{GR} that for each $\epsilon>0$, by choosing $\delta$ small and $N$ large enough, the probability that the walkers move distance $N^\kappa L_N$ apart 
without getting within distance $L_N$ from each other is at least $c\delta N^{-\kappa-\epsilon}$.

Hence there is a polynomially small probability of making an
excursion of size $N^\kappa L_N$ before returning to a distance
$L_N$. On the other hand once we have such a big excursion
Lemma \ref{lem:largedev} implies that it will take more than 
$N$ steps to come back, indeed 
\begin{displaymath}
  \begin{split}
    \bP^2&\left(\left\{\inf_{\ell+1\leq j\leq \ell+N} \|X_j-Y_j\|\leq L_N\right\}\;\bigg|
      \;\|X_\ell-Y_\ell\|\geq N^\kappa L_N\right)\\
    &  \leq CNe^{-CN^{2\kappa-1}}.
  \end{split}
\end{displaymath}
The last two estimates imply Lemma \ref{lem:Ex}
\end{proof}

\subsection{Comparison Lemma}\label{sec:comparison}
\begin{proof}[\bfseries Proof of Lemma \ref{PrComp}]
Let $U_1, U_2\dots  U_n\dots $ be random variables which are independent and uniformly 
distributed on $[0,1].$ Define
$$ \xi_n^*=-1 \text{ if } U_n<P(\xi_n=-1|\xi_1=\xi^*_1, \dots, \xi_{n-1}=\xi_{n-1}^*)
\text{ and } \xi_n=1 \text{ otherwise} . $$
Also let
$\txi_n^*=-1$ if $U_n<1-p$ and $\txi_n^*=1$ otherwise. Let
$$ \cX_n^*=\sum_{j=1}^n \xi_j^*, \quad \tcX_n^*=\sum_{j=1}^n \txi_j^*. $$
Then $\{\cX^*_n\}$ has the same distribution as $\{\cX_n\},$
$\{\tcX^*_n\}$ has the same distribution as $\{\tcX_n\}$ and
$\cX^*_n\geq \tcX^*_n.$
\end{proof}

\end{document}